\documentclass[amstex,16pt, amssymb]{article}
\usepackage{mathtext}
\usepackage[cp1251]{inputenc}
\usepackage[T2A]{fontenc}
\usepackage[dvips]{graphicx}
\usepackage{amsmath}
\usepackage{amssymb}
\usepackage{amsxtra}
\usepackage{latexsym}
\usepackage{ifthen}
\usepackage{amsthm}

\def\Xint#1{\mathchoice
   {\XXint\displaystyle\textstyle{#1}}
   {\XXint\textstyle\scriptstyle{#1}}
   {\XXint\scriptstyle\scriptscriptstyle{#1}}
   {\XXint\scriptscriptstyle\scriptscriptstyle{#1}}
   \!\int}
\def\XXint#1#2#3{{\setbox0=\hbox{$#1{#2#3}{\int}$}
     \vcenter{\hbox{$#2#3$}}\kern-.5\wd0}}
\def\dashint{\Xint-}

\numberwithin{equation}{section}

\begin{document}

\title{On Dirichlet problem for degenerate \\ Beltrami equations with sources}

\author{V. Gutlyanski\u\i{}, O. Nesmelova, V. Ryazanov, E. Yakubov}







\date{}

\maketitle

\begin{abstract}
The present paper is devoted to the study of the Dirichlet problem
${\rm{Re}}\,\omega(z)\to\varphi(\zeta)$ as $z\to\zeta,$ $z\in
D,\zeta\in \partial D,$ with continuous boundary data $\varphi
:\partial D\to\mathbb R$ for Beltrami equations
$\omega_{\bar{z}}=\mu(z) \omega_z+\sigma (z)$, $|\mu(z)|<1$ a.e.,
with sources $\sigma :D\to\mathbb C$ in the case of locally uniform
ellipticity. In this case, we establish a series of effective
integral criteria of the type of BMO, FMO, Calderon-Zygmund, Lehto
and Orlicz on singularities of the equations at the boundary for
existence, representation and regularity of its solutions in
arbitrary bounded domains $D$ of the complex plane $\mathbb C$ with
no boun\-da\-ry component degenerated to a single point for sources
$\sigma$ in $L_p(D)$, $p>2$, with compact support in $D$. Moreover,
we prove in such domains existence, representation and regularity of
weak solutions of the Dirichlet problem for the Poisson type
equation ${\rm div} [A(z)\nabla\,u(z)] = g(z)$ whose source $g\in
L_p(D)$, $p>1$, has compact support in $D$ and whose mat\-rix valued
coefficient $A(z)$ guarantees its locally uniform ellipticity.
\end{abstract}

\bigskip


\bigskip



\section{Introduction}

Let $D$ be a domain in the complex plane $\mathbb C.$ We study the
{\bf Dirichlet problem}
\begin{equation}\label{eqDIRICHLET}\lim\limits_{z\to\zeta}{\rm
Re}\,\omega(z)=\varphi(\zeta)\qquad\forall\ \zeta\in\partial D\
,\end{equation} see e.g. \cite{Bojar} and \cite{Vek}, with
continuous boundary data $\varphi :\partial D\to\mathbb R$  in
arbitrary bounded domains $D$ with no boun\-da\-ry component
degenerated to a single point for the {\bf inhomogeneous Beltrami
equation}
\begin{equation}\label{eqBeltrami}
\omega_{\bar{z}}\ =\ \mu(z)\cdot \omega_z\ +\ \sigma(z)\ ,\ \ \ \
z\in D\ ,
\end{equation}
{\bf with a source} $\sigma :D\to\mathbb C$ in $L_p$, $p>2$, where
$\mu : D\to\mathbb C$ is a  measurable function with $|\mu(z)|< 1$
a.e., $\omega_{\bar z}=(\omega_x+i\omega_y)/2$,
$\omega_{z}=(\omega_x-i\omega_y)/2$, $z=x+iy$, $\omega_x$ and
$\omega_y$ are partial derivatives of the function $\omega$ in $x$
and $y$, respectively.

Moreover, in general bounded domains $D$ without any boundary
component degenerated to a single point, we prove the
cor\-res\-pon\-ding theorems on existence, representation and
re\-gu\-la\-ri\-ty of solutions for the {\bf classical Dirichlet
problem}
\begin{equation}\label{eqDIR}
\lim_{z\to\zeta}\ u(z)\ =\ \varphi(\zeta)\ \ \ \ \ \forall\
\zeta\in\partial D
\end{equation}
with continuous boundary data $\varphi :\partial D\to\mathbb R$ to
the {\bf Poisson type equation}
\begin{equation}\label{eqPotential}
{\rm div} [A(z)\nabla\,u(z)] = g(z)
\end{equation}
with a source $g:D\to\mathbb R$ in $L_p(D)$, $p>1$, see e.g.
background for $g\equiv 0$ in \cite{GRSY*}.

The request on domains to have no boundary component degenerated to
a single point is necessary. Indeed, consider the punctured unit
disk $\mathbb D_0:=\mathbb D\setminus\{ 0\}$. Setting
$\varphi(\zeta)\equiv 1$ on $\partial\mathbb D$ and $\varphi(0)= 0$,
we see that $\varphi$ is continuous on $\partial\mathbb
D_0=\partial\mathbb D\cup\{ 0\}$. Let us assume that there is a
harmonic function $u$ satisfying (\ref{eqDIR}). Then $u$ is bounded
by the maximum principle for harmonic functions and by the classical
Cauchy--Riemann theorem, see also Theorem V.4.2 in \cite{Ne}, the
extended $u$ is harmonic in $\mathbb D$. Thus, by
Mean-Value-Property for harmonic functions we disprove the above
assumption, see e.g. Theorem 0.2.4 in \cite{ST}.

In this connection, recall that a boundary point $p$ of a domain $D$
in $\mathbb R^n, n\geq 2,$ is called regular if each solution of the
Dirichlet problem for the Laplace equation in $D$, whose Dirichlet
boundary date is continuous at $p$, is also continuous at $p$. The
famous Wiener criterion for regularity of a boun\-da\-ry point, see
\cite{Wien}, that has been formulated in terms of so-called barrier
functions, generally speaking, has  no satisfactory geometric
interpretation. However, there is a very simple geometric criterion
of regular points in the case of $\mathbb C$. Namely, a point
$p\in\partial D$ is regular if $p$ belongs to a component of
$\partial D$ that is not degenerated to a single point, see Theorem
4.2.2 in \cite{Rans}. The above example shows that this condition is
not only sufficient but also necessary for regularity of a boundary
point. Thus, results on the Dirichlet problem (\ref{eqDIR}) to the
Poisson type equations (\ref{eqPotential}) given in the end of the
paper were obtained for the most general ad\-mis\-sib\-le domains.

For the case $\|\mu\|_{\infty}<1$, (\ref{eqBeltrami}) was first
introduced by L. Ahlfors and L. Bers in the paper \cite{ABe}, see
also the Ahlfors monograph \cite{Alf}. Here we study the case of
locally uniform ellipticity of the equation (\ref{eqBeltrami}) when
its {\bf dilatation quotient} $K_{\mu}$ is bounded only locally in
$D$,
\begin{equation}\label{K}
K_{\mu}(z)\ :=\ \frac{1\ +\ |\mu(z)|}{1\ -\ |\mu(z)|}\ ,
\end{equation}
i.e., if $K_{\mu}\in L_{\infty}$ on each compact set in $D$ but
admits singularities at the boundary. A point $\zeta\in\partial D$
is called a {\bf singular point of the equation} (\ref{eqBeltrami})
if $K_{\mu}\notin L_{\infty}$ on each neighborhood of the point.

First of all, we show that, if $D$ is an arbitrary simply connected
domain in $\mathbb C$, then the Dirichlet problem
(\ref{eqDIRICHLET}) to the equation (\ref{eqBeltrami}) has a
solution $\omega$ in class $W^{1,2}_{\rm loc}(D)$ for a wide circle
of singularities of (\ref{eqBeltrami}) at the boundary. Moreover, it
is unique up to an additive pure imaginary constant, and it can be
represented through the so-called generalized analytic functions
with sources.

Recall that the Vekua monograph \cite{Vek} was devoted to {\bf
generalized analytic functions}, i.e., continuous complex valued
functions $H(z)$ of one complex variable $z=x+iy$ of class
$W^{1,1}_{\rm loc}$ in a domain $D$ satisfying the equations
\begin{equation}\label{eqG}
\partial_{\bar z}H\ +\ aH\ +\ b {\overline H}\ =\ S\ ,\ \ \  \partial_{\bar
z}\ :=(\partial_x+i\partial_y)/2
\end{equation}
with complex valued coefficients  $a,b,S\in L_{p}(D)$, $p>2$.

\medskip

The papers \cite{GNRY} and \cite{R4} were devoted to boundary value
problems with measurable data for the spacial case of {\bf
generalized analytic functions H with sources} $S:D\to\mathbb C$,
when $a\equiv 0\equiv b$,
\begin{equation}\label{eqS}
\partial_{\bar z}H(z)\ =\ S(z)\ ,\ \ \ z\in D\ ,
\end{equation}
in regular enough domains $D$.


Then we establish here similar theorems on existence, representation
and re\-gu\-la\-ri\-ty of the so-called multi-valued solutions of
the Dirichlet problem (\ref{eqDIRICHLET}) with continuous boundary
data $\varphi:\partial D\to\mathbb R$ to the Beltrami equations with
sources (\ref{eqBeltrami}) in arbitrary domains $D$ in $\mathbb C$
without any boundary components degenerated to a single point.
Finally, we resolve the classical Dirichlet problem (\ref{eqDIR}) to
Poisson type equations (\ref{eqPotential}) in the general domains.


In particular, we give here a representation of the given solutions
for (\ref{eqDIR}) to (\ref{eqPotential}) through generalized
harmonic functions with  sources. Recall that the paper \cite{R4}
were devoted to the existence of nonclassical continuous solutions
of class $W^{2,p}_{\rm loc}$ to various boundary-value problems with
arbitrary boundary data that were measurable with respect to the
length measure in domains with rectifiable boundaries for {\bf
generalized harmonic functions with sources} $G:D\to\mathbb R$ in
$L_p(D)$, $p>2$, satisfying the Poisson equations
\begin{equation}\label{eqSSS}
\triangle U(z)\ =\ G(z)\ .
\end{equation} Note that by the Sobolev embedding theorem, see
Theorem I.10.2 in \cite{So}, such functions $U$ belong to the class
$C^1$. Similar results were also proved in \cite{GNRY} with
arbitrary boundary data that were measurable with respect to the
logarithmic capacity in special domains with non\-rec\-ti\-fi\-able
boundaries. In the case $G$ in $L_p(D)$, $p>1$, we will call
continuous solutions of the Poisson equation (\ref{eqSSS}) of the
class $W^{2,p}_{\rm loc}$ {\bf weak generalized harmonic functions
with the sources} $G$.

\medskip

The paper is organized as follows. Sections 2 and 3  contain the
main lemma and a series of other effective integral criteria,
respectively, for existence, representation and regularity of
solutions of the Dirichlet problem (\ref{eqDIRICHLET}) to the
Beltrami equations with sources (\ref{eqBeltrami}) in the case of
arbitrary simply connected domains. Sections 4 and 5 include similar
results on multi-valued solutions of the Dirichlet problem
(\ref{eqDIRICHLET}) to the equations (\ref{eqBeltrami}) in the case
of domains $D$ with no boundary components degenerated to a single
point. Finally, in arbitrary such domains we obtain existence,
representation and regularity results for the classsical Dirichlet
problem (\ref{eqDIR}) to Poisson type equations (\ref{eqPotential})
in Section 6.

\section{The main lemma in simply connected domains}

It is well known that the homogeneous Beltrami equation
\begin{equation}\label{a}
f_{\bar{z}}=\mu(z) f_z
\end{equation}
is the basic equation in analytic theory of quasiconformal and
quasiregular mappings in the plane with numerous applications in
nonlinear elasticity, gas flow, hydrodynamics  and other sections of
natural sciences. Note that continuous functions $f$ with
generalized derivative by Sobolev $f_{\bar z}=0$ are analytic
functions, see e.g. Lemma 1 in \cite{ABe}, that corresponds to the
case $\mu(z)\equiv 0$ in (\ref{a}).

The equation (\ref{a}) is called {\bf degenerate} if ${\rm
ess}\,{\rm sup}\,K_{\mu}(z)=\infty$. It is known that if $K_{\mu}$
is bounded, then the equation has ho\-meo\-mor\-phic solutions in $
W^{1,2}_{\rm loc}$, see e.g. monographs \cite{Alf}, \cite{BGMR} and
\cite{LV}, called {\bf quasiconformal mappings}. Recently, a series
of effective criteria for existence of homeomorphic solutions in $
W^{1,1}_{\rm loc}$ have been also established for the degenerate
Beltrami equations, see e.g. historic comments with relevant
references in monographs \cite{AIM}, \cite{GRSY} and \cite{MRSY}.

These criteria were formulated both in terms of $K_{\mu}$ and the
more refined quantity that takes into account not only the modulus
of the complex coefficient $\mu$ but also its argument
\begin{equation}\label{eqTangent} K^T_{\mu}(z,z_0)\ :=\
\frac{\left|1-\frac{\overline{z-z_0}}{z-z_0}\mu (z)\right|^2}{1-|\mu
(z)|^2} \end{equation} called {\bf tangent dilatation quotient} of
Beltrami equations with respect to a point $z_0\in\mathbb C$, see
e.g. \cite{And}, \cite{BGR$_1$}, \cite{BGR$_2$}, \cite{GMSV},
\cite{GRSY}, \cite{Le}, \cite{MRSY} and \cite{RSY}-\cite{RSY$_5$}.
Note that
\begin{equation}\label{eqConnect} K^{-1}_{\mu}(z)\leqslant K^T_{\mu}(z,z_0) \leqslant K_{\mu}(z)
\ \ \ \ \ \ \ \forall\ z\in D\, ,\ z_0\in \Bbb C\ .\end{equation}

Let $D$ be a domain in the complex plane ${\Bbb C}$. A function
$f:D\to\mathbb C$ in the Sobolev class $W^{1,1}_{\rm loc}$ is called
a {\bf regular solution} of the Beltrami equation (\ref{a}) if $f$
satisfies (\ref{a}) a.e. and its Jacobian $J_f(z)=|f_z|^2-|f_{\bar
z}|^2>0$ a.e. By Lemma 3 and Remark 2 in \cite{RSY$_4$}, we have the
following statement on the existence of regular homeomorphic
solutions $f$ in $\mathbb C$ for
the Beltrami equation (\ref{a}).

\medskip

{\bf Proposition 1.} {\it Let $\mu : \mathbb C\to{\Bbb C}$ be a
measurable function with $|\mu (z)| < 1$  a.e. and $K_{\mu}\in
L_{1,\rm loc}(\mathbb C).$ Suppose that, for each $z_0\in \mathbb C$
with some $\varepsilon_0=\varepsilon(z_0)>0,$
\begin{equation}\label{eqT}
\int\limits_{{\varepsilon}<|z-z_0|<{\varepsilon}_0}\
K^T_{{\mu}}(z,z_0)\cdot{\psi}^2_{z_0,{\varepsilon}}(|z-z_0|)\ dm(z)\
=\ o(I^2_{z_0}({\varepsilon})) \ \ \ \hbox{as ${\varepsilon}\to 0$}
\end{equation}
for a family of measurable functions
${\psi}_{z_0,{\varepsilon}}:(0,\varepsilon_0)\to(0,\infty),$
${\varepsilon}\in(0,{\varepsilon}_0),$ with
\begin{equation}\label{eqI}
I_{z_0}({\varepsilon})\ \colon =\
\int\limits_{{\varepsilon}}^{{\varepsilon}_0}{\psi}_{z_0,{\varepsilon}}(t)\
dt\ <\ \infty\ \ \ \ \ \ \forall\
{\varepsilon}\in(0,{\varepsilon}_0)\ .\end{equation} Then the
Beltrami equation (\ref{a}) has a regular homeomorphic solution
$f^{\mu}$.}

\medskip

Here $dm(z)$ corresponds to the Lebesgue measure in ${\Bbb C}$ and
by (\ref{eqConnect}) $K^T_{\mu}$ can be replaced by $K_{\mu}$. We
call such solutions $f^{\mu}$ of (\ref{a}) {\bf $\mu -$conformal
mappings}.

\medskip

{\bf Lemma 1.} {\it Let $D$ be a bounded simply connected domain in
${\Bbb C}$, $\sigma\in L_p(D)$, $p>2$, with compact support in $D$,
$\mu:D\to{\Bbb C}$ be a measurable function with $|\mu(z)|<1$ a.e.,
$K_{\mu}$ be locally bounded in $D$, $K_{\mu}\in L_1(D)$ and
conditions (\ref{eqT}) and (\ref{eqI}) hold for all
$z_0\in\partial{D}$.

Then the Beltrami equation (\ref{eqBeltrami}) with the source
$\sigma$ has a locally H\"older continuous solution $\omega$ in the
class $W^{1,2}_{\rm loc}$ of the Dirichlet problem
(\ref{eqDIRICHLET}) in $D$ for each continuous function
$\varphi:\partial D\to{\Bbb R}$ that is unique up to an additive
pure imaginary constant.

Moreover, $\omega =h\circ f$, where $f:\mathbb C\to\mathbb C$ is a
$\mu -$conformal mapping with $\mu$ extended by zero outside of $D$
and $h:D_*\to\mathbb C$  is a generalized analytic function in
$D_*:=f(D)$ with the source $S$ of the class $L_{p_*}(D_*)$ for some
$p_*\in(2,p)$,
\begin{equation}\label{r}
S\ :=\ \sigma\cdot\frac{f_z}{J}\circ f^{-1}\ ,
\end{equation} where $J=|f_z|^2-|f_{\bar z}|^2$ is the Jacobian of $f$, that satisfies the Dirichlet condition
\begin{equation}\label{eqHOLOMORPHIC}
\lim_{w\to\zeta}\ {\rm Re}\, h(w)\ =\ \varphi_*(\zeta)\ \ \ \ \
\forall\ \zeta\in\partial D_*\ ,\ \ \ \ \ \mbox{where
$\varphi_*:=\varphi\circ f^{-1}|_{\partial D_*}$.}
\end{equation}}

\medskip

We assume here and further that the dilatation quotients
$K^T_{\mu}(z,z_0)$ and $K_{\mu}(z)$ are extended by $1$ outside of
the domain $D$.

\medskip

{\bf Remark 1.} In tern, the generalized analytic function $h$ with
the source $S$ by Theorem 1.16 in \cite{Vek} has the representation
$h = {\cal A} + H$, where
\begin{equation}\label{eqINTEGRAL}
H(w)\ =\ -\frac{1}{\pi}\int\limits_{D_*}\frac{S(\zeta)}{\zeta -w}\,
dm(\zeta)\ ,\ \ \ w\in\mathbb C\ ,
\end{equation}
with $H_{\overline w}=S$, and ${\cal A}$ is a holomorphic function
in $D_*$ with the Dirichlet condition
\begin{equation}\label{eqSTAR}
\lim_{w\to\zeta}\ {\rm Re}\, {\cal A}(w)\ =\ \varphi^*(\zeta)\ \ \ \
\ \forall\ \zeta\in\partial D_*\ ,\ \ \  \mbox{where
$\varphi^*:=\varphi_* - {\rm Re}\, H|_{\partial D_*}$\ .}
\end{equation}
Note that $H$ is $\alpha_* -$H\"older continuous in $D_*$ with
$\alpha_* =1-2/p_*$ by Theorem 1.19 and $H|_{D_*}\in W^{1,p_*}(D_*)$
by Theorems 1.36 and 1.37 in \cite{Vek}. Note also that $f$ and
$f^{-1}$ are locally quasiconformal mappings in $D$ and $D_*$,
respectively.


\begin{proof}
Let us first show the uniqueness of the desired solution. Indeed, if
$\omega_1$ and $\omega_2$ are such solutions, then $\Omega
:=\omega_2-\omega_1$ is such a solution for the Dirichlet problem
with zero boundary date to the homogeneous Beltrami equation
(\ref{a}). Consider the function $\Omega^*:=\Omega\circ f^{-1}$.
First of all, note that $f^*:=f^{-1}|_{D_*}$ is a locally
quasiconformal mapping and, in particular, $f^*\in W^{1,2}_{\rm
loc}(D_*)$. Hence
$$
\Omega^*_{\overline w}\ =\ \Omega_{z}\circ f^*\cdot f^*_{\overline
w}\ +\ \Omega_{\overline{z}}\circ f^*\cdot\overline{f^*_{w}}\ =\
\Omega_{z}\circ f^*\cdot\left[ \ f^*_{\overline w}\ +\ \mu\circ
f^*\cdot\overline{f^*_{w}}\ \right]\ \ \ \hbox{a.e. in $D_*$}
$$
and $\Omega^*\in W^{1,1}_{\rm loc}$ by Lemma III.6.4 in \cite{LV}.
Consequently, $\Omega^*_{\overline w}=0$ a.e. in $D_*$ by I.C(3) in
\cite{Alf}. Thus, the function $\Omega^*$ is analytic in $D_*$ by
Lemma 1 in \cite{ABe}. Its real part $u$ satisfies zero Dirichlet
condition and by the maximum principle for harmonic functions
$u\equiv 0$ in $D_*$. Thus, $\Omega$ is a pure imaginary constant.

Now, let us prove that in (\ref{r}) $S\in L_{p_*}(D_*)$ for some
$p_*\in (2,p)$. Indeed, let ${\tilde D}$ be a subdomain of $D$ with
its closure in $D$, containing the support of $\sigma$. Setting
${\tilde \mu}\equiv\mu$ on ${\tilde D}$ and zero for the rest points
of $\mathbb C$, we see that $K_{\tilde\mu}\in L_{\infty}(\mathbb
C)$. Consequently, we obtain that the function ${\tilde S}:{\tilde
D}_*\to\mathbb C$, ${\tilde D}_*:=f^{\tilde\mu}(D_*)$,
\begin{equation}\label{R}
{\tilde S}\ :=\
\left(\frac{f^{\tilde\mu}_z}{J^{\tilde\mu}}\cdot\sigma\right)\circ
\left(f^{\tilde\mu}\right)^{-1}\ ,
\end{equation} where $f^{\tilde\mu}:\mathbb C\to\mathbb C$ is the
${\tilde\mu} -$ conformal mapping from Theorem B in \cite{GNRY*} and
$J^{\tilde\mu}$ is the Jacobian of $f^{\tilde\mu}$, belongs to class
$L_{p_*}({\tilde D}_*)$ for some $p_*\in (2,p)$ by Remark 2 to Lemma
1 in \cite{GNRY*}. However, $f^{\tilde\mu}|_{\tilde D}=C\circ
f|_{\tilde D}$, where $C$ is a conformal mapping on $
D^{\prime}_*:=f(\tilde D)$ because both $f$ and $f^{\tilde\mu}$ are
quasiconformal mappings on $\tilde D$ with the same complex
characteristic $\mu$ there. Consequently,
\begin{equation}
S\ =\ \left({\tilde S}\cdot \overline{C^{\prime}}\right)\circ C
\end{equation}
on $D^{\prime}_*\subset D_*$ containing the compact support of the
function $\Sigma := \sigma\circ f^{-1}$. Thus, really $S\in
L_{p_*}(D_*)$ for the same $p_*\in (2,p)$.

Next, let $\varphi^*:=\left(\varphi\circ f^{-1} - {\rm Re}\,
H\right)|_{\partial D_*}$, where $H$ is the generalized analytic
function (\ref{eqINTEGRAL}) with the source $S$. Then by Corollary
4.1.8 and Theorem 4.2.1 in \cite{Rans} there is a harmonic function
$u:D_*\to\mathbb R$ satisfying the Dirichlet condition
\begin{equation}\label{eqHARMONIC}
\lim_{w\to\zeta}\ u(w)\ =\ \varphi^*(\zeta)\ \ \ \ \ \forall\
\zeta\in\partial D_*
\end{equation}
because $D_*$ is a bounded simple connected domain that, of course,
has at least 2 boundary points. On the other hand, there is its
conjugate harmonic function $v:D_*\to\mathbb R$ such that ${\cal
A}:=u+iv:D_*\to\mathbb C$ forms a holomorphic function again because
of the domain $D_*$ is simply connected, see e.g. arguments in the
beginning of the book \cite{Ko}. Thus, the function $\omega:=h\circ
f$, where $h:={\cal A}+H$, gives the desired solution of the
Dirichlet problem (\ref{eqDIRICHLET}) in $D$ to the Beltrami
equation (\ref{eqBeltrami}) with the source $\sigma$ of the class
$W^{1,1}_{\rm loc}$ by Lemma III.6.4 in \cite{LV}, see also Remark
1. Arguing as in the previous item with the application of the
auxiliary quasiconformal mapping $f^{\tilde\mu}:\mathbb C\to\mathbb
C$, it is easy to prove on the basis of Lemma 1 in \cite{GNRY*} that
$\omega\in W^{1,2}_{\rm loc}$. Finally, the given solution $\omega$
is locally H\"older continuous because the function $h$ and the
mapping $f$ are so, see Remark 1.
\end{proof}


{\bf Remark 2.} Note that if the family of the functions
$\psi_{z_0,\varepsilon}(t)\equiv\psi_{z_0}(t)$, $z_0\in\partial{D}$,
in Lemma 1 is independent on the parameter $\varepsilon$, then the
condition (\ref{eqT}) implies that $I_{z_0}(\varepsilon)\to \infty$
as $\varepsilon\to 0$. This follows immediately from arguments by
contradiction, apply for it (\ref{eqConnect}) and the condition
$K_{\mu}\in L_1(D)$. Note also that (\ref{eqT}) holds, in
particular, if, for some $\varepsilon_0=\varepsilon(z_0)$,
\begin{equation}\label{333omal}
\int\limits_{|z-z_0|<\varepsilon_0}
K^T_{\mu}(z,z_0)\cdot\psi_{z_0}^2(|z-z_0|)\,dm(z)<\infty \qquad
\forall\ z_0\in\partial{D}\end{equation} and
$I_{z_0}(\varepsilon)\to \infty$ as $\varepsilon\to 0$. In other
words, for the solvability of the Dirichlet problem
(\ref{eqDIRICHLET}) in $D$ for the Beltrami equations with sources
(\ref{eqBeltrami}) for all con\-ti\-nu\-ous boundary functions
$\varphi$, it is sufficient that the integral in (\ref{333omal})
converges for some nonnegative function $\psi_{z_0}(t)$ that is
locally integrable over $(0,\varepsilon_0 ]$ but has a nonintegrable
singularity at $0$. The functions $\log^{\lambda}(e/|z-z_0|)$,
$\lambda\in (0,1)$, $z\in\Bbb D$, $z_0\in\overline{\Bbb D}$, and
$\psi(t)=1/(t \,\, \log(e/t))$, $t\in(0,1)$, show that the condition
(\ref{333omal}) is compatible with the condition
$I_{z_0}(\varepsilon)\to\infty $ as $\varepsilon\to 0$. Furthermore,
the condition (\ref{eqT}) shows that it is sufficient for the
solvability of the Dirichlet problem even if the integral in
(\ref{333omal}) is divergent in a controlled way.

\section{The main criteria in simply connected domains}

Lemma 1 makes it to be possible to derive a number of effective
integral criteria for solvability of the Dirichlet problem to the
Beltrami equations with sources.

\bigskip

Recall first that a real-valued function $u$ in a domain $D$ in
${\Bbb C}$ is said to be of {\bf bounded mean oscillation} in $D$,
abbr. $u\in{\rm BMO}(D)$, if $u\in L_{\rm loc}^1(D)$ and
\begin{equation}\label{lasibm_2.2_1}\Vert u\Vert_{*}:=
\sup\limits_{B}{\frac{1}{|B|}}\int\limits_{B}|u(z)-u_{B}|\,dm(z)<\infty\,,\end{equation}
where the supremum is taken over all discs $B$ in $D$ and
$$u_{B}={\frac{1}{|B|}}\int\limits_{B}u(z)\,dm(z)\,.$$ We write $u\in{\rm BMO}_{\rm loc}(D)$ if
$u\in{\rm BMO}(U)$ for each relatively compact subdomain $U$ of $D$.
We also write sometimes for short BMO and ${\rm BMO}_{\rm loc }$,
respectively.

The class BMO was introduced by John and Nirenberg (1961) in the
paper \cite{JN} and soon became an important concept in harmonic
analysis, partial differential equations and related areas, see e.g.
\cite{HKM} and \cite{RR}.

\medskip

A function $\varphi$ in BMO is said to have {\bf vanishing mean
oscillation}, abbr. $\varphi\in{\rm VMO}$, if the supremum in
(\ref{lasibm_2.2_1}) taken over all balls $B$ in $D$ with
$|B|<\varepsilon$ converges to $0$ as $\varepsilon\to0$. VMO has
been introduced by Sarason in \cite{Sar}. There are a number of
papers devoted to the study of partial differential equations with
coefficients of the class VMO, see e.g. \cite{CFL}, \cite{IS},
\cite{MRV}, \cite{Pal}, \cite{Ra$_1$} and \cite{Ra$_2$}.

\medskip

{\bf Remark 3.} Note that $W^{\,1,2}\left({{D}}\right) \subset VMO
\left({{D}}\right),$ see e.g. \cite{BN}.

\medskip

Following \cite{IR}, we say that a function $\varphi:D\to{\Bbb R}$
has {\bf finite mean oscillation} at a point $z_0\in D$, abbr.
$\varphi\in{\rm FMO}(z_0)$, if
\begin{equation}\label{FMO_eq2.4}\overline{\lim\limits_{\varepsilon\to0}}\ \ \
\dashint_{B(z_0,\varepsilon)}|{\varphi}(z)-\widetilde{\varphi}_{\varepsilon}(z_0)|\,dm(z)<\infty\,,\end{equation}
where \begin{equation}\label{FMO_eq2.5}
\widetilde{\varphi}_{\varepsilon}(z_0)=\dashint_{B(z_0,\varepsilon)}
{\varphi}(z)\,dm(z)\end{equation} is the mean value of the function
${\varphi}(z)$ over the disk $B(z_0,\varepsilon):=\{ z\in\mathbb C:
|z-z_0|<\varepsilon\}$. Note that the condition (\ref{FMO_eq2.4})
includes the assumption that $\varphi$ is integrable in some
neighborhood of the point $z_0$. We say also that a function
$\varphi:D\to{\Bbb R}$ is of {\bf finite mean oscillation in $D$},
abbr. $\varphi\in{\rm FMO}(D)$ or simply $\varphi\in{\rm FMO}$, if
$\varphi\in{\rm FMO}(z_0)$ for all points $z_0\in D$. We write
$\varphi\in{\rm FMO}(\overline{D})$ if $\varphi$ is given in a
domain $G$ in $\Bbb{C}$ such that $\overline{D}\subset G$ and
$\varphi\in{\rm FMO}(G)$.

\medskip

The following statement is obvious by the triangle inequality.

\medskip

{\bf Proposition 2.} {\it If, for a  collection of numbers
$\varphi_{\varepsilon}\in{\Bbb R}$,
$\varepsilon\in(0,\varepsilon_0]$,
\begin{equation}\label{FMO_eq2.7}\overline{\lim\limits_{\varepsilon\to0}}\ \ \
\dashint_{B(z_0,\varepsilon)}|\varphi(z)-\varphi_{\varepsilon}|\,dm(z)<\infty\,,\end{equation}
then $\varphi $ is of finite mean oscillation at $z_0$.}

\medskip

In particular, choosing here  $\varphi_{\varepsilon}\equiv0$,
$\varepsilon\in(0,\varepsilon_0]$ in Proposition 1, we obtain the
following.

\medskip

{\bf Corollary 1.} {\it If, for a point $z_0\in D$,
\begin{equation}\label{FMO_eq2.8}\overline{\lim\limits_{\varepsilon\to 0}}\ \ \
\dashint_{B(z_0,\varepsilon)}|\varphi(z)|\,dm(z)<\infty\,,
\end{equation} then $\varphi$ has finite mean oscillation at
$z_0$.}

\medskip

Recall that a point $z_0\in D$ is called a {\bf Lebesgue point} of a
function $\varphi:D\to{\Bbb R}$ if $\varphi$ is integrable in a
neighborhood of $z_0$ and \begin{equation}\label{FMO_eq2.7a}
\lim\limits_{\varepsilon\to 0}\ \ \ \dashint_{B(z_0,\varepsilon)}
|\varphi(z)-\varphi(z_0)|\,dm(z)=0\,.\end{equation} It is known
that, almost every point in $D$ is a Lebesgue point for every
function $\varphi\in L_1(D)$. Thus, we have by Proposition 1 the
next corollary.

\medskip

{\bf Corollary 2.} {\it Every locally integrable function
$\varphi:D\to{\Bbb R}$ has a finite mean oscillation at almost every
point in $D$.}

\medskip

{\bf Remark 4.} Note that the function
$\varphi(z)=\log\left(1/|z|\right)$ belongs to BMO in the unit disk
$\Delta$, see, e.g., \cite{RR}, p. 5, and hence also to FMO.
However, $\widetilde{\varphi}_{\varepsilon}(0)\to\infty$ as
$\varepsilon\to0$, showing that condition (\ref{FMO_eq2.8}) is only
sufficient but not necessary for a function $\varphi$ to be of
finite mean oscillation at $z_0$. Clearly, ${\rm BMO}(D)\subset{\rm
BMO}_{\rm loc}(D)\subset{\rm FMO}(D)$ and as well-known ${\rm
BMO}_{\rm loc}\subset L_{\rm loc}^p$ for all $p\in[1,\infty)$, see,
e.g., \cite{JN} or \cite{RR}. However, FMO is not a subclass of
$L_{\rm loc}^p$ for any $p>1$ but only of $L_{\rm loc}^1$. Thus, the
class FMO is much more wider than ${\rm BMO}_{\rm loc}$.

\medskip

Versions of the next lemma has been first proved for the class BMO
in \cite{RSY}. For the FMO case, see the papers \cite{IR, RS,
RSY$_2$, RSY$_3$}  and the monographs \cite{GRSY} and \cite{MRSY}.

\medskip

{\bf Lemma 2.} {\it Let $D$ be a domain in ${\Bbb C}$ and let
$\varphi:D\to{\Bbb R}$ be a  non-negative function  of the class
${\rm FMO}(z_0)$ for some $z_0\in D$. Then
\begin{equation}\label{eq13.4.5}\int\limits_{\varepsilon<|z-z_0|<\varepsilon_0}\frac{\varphi(z)\,dm(z)}
{\left(|z-z_0|\log\frac{1}{|z-z_0|}\right)^2}=O\left(\log\log\frac{1}{\varepsilon}\right)\
\quad\text{as}\quad\varepsilon\to 0\end{equation} for some
$\varepsilon_0\in(0,\delta_0)$ where $\delta_0=\min(e^{-e},d_0)$,
$d_0=\sup\limits_{z\in D}|z-z_0|$.}


Recall that we assume further that the dilatation quotients
$K^T_{\mu}(z,z_0)$ and $K_{\mu}(z)$ are extended by $1$ outside of
the domain $D$.

\medskip

Choosing $\psi(t)=1/\left(t\, \log\left(1/t\right)\right)$ in Lemma
1, see also Remark 1, we obtain by Lemma 2 the following result with
the FMO type criterion.

\medskip

{\bf Theorem 1.} {\it Let $D$ be a bounded simply connected domain
in $\mathbb C$, $\sigma\in L_p(D)$, $p>2$, with compact support,
$\mu:D\to{\Bbb C}$ be a measurable function with $|\mu(z)|<1$ a.e.,
$K_{\mu}$ be locally bounded in $D$, $K_{\mu}\in L_1(D)$,
$K^T_{\mu}(z,z_0)\leqslant Q_{z_0}(z)$ a.e. in $U_{z_0}$ for each
point $z_0\in \partial{D}$, a neighborhood $U_{z_0}$ of $z_0$, a
function $Q_{z_0}: U_{z_0}\to[0,\infty]$ in the class ${\rm
FMO}({z_0})$.

Then the Beltrami equation (\ref{eqBeltrami}) with the source
$\sigma$ has a locally H\"older continuous solution $\omega$ in the
class $W^{1,2}_{\rm loc}$ of the Dirichlet problem
(\ref{eqDIRICHLET}) in $D$ for each continuous function
$\varphi:\partial D\to{\Bbb R}$ that is unique up to an additive
pure imaginary constant.

Moreover, $\omega =h\circ f$, $h := {\cal A} + H$, where $f:\mathbb
C\to\mathbb C$ is a $\mu -$conformal mapping with $\mu$ extended by
zero outside of $D$, $H:D_*\to\mathbb C$  is a generalized analytic
function in $D_*:=f(D)$ with the source $S$ calculated in (\ref{r})
and ${\cal A}$ is a holomorphic function in $D_*$ with the Dirichlet
condition (\ref{eqSTAR}).}

\medskip

By Corollary 1 we obtain the following nice consequence of Theorem
1, where $B(z_0,\varepsilon)$ denote the infinitesimal disks $\{ z
\in\mathbb C:\, |z-z_0|\,<\, \varepsilon\}$ centered at
$z_0\in\partial D$.

\medskip

{\bf Corollary 3.} {\it Let $D$ be a bounded simply connected domain
in $\mathbb C$, $\sigma\in L_p(D)$, $p>2$, with compact support,
$\mu:D\to{\Bbb C}$ be a measurable function with $|\mu(z)|<1$ a.e.,
$K_{\mu}$ be locally bounded in $D$, $K_{\mu}\in L_1(D)$ and, for
each point $z_0\in\partial{D}$,
\begin{equation}\label{eqMEAN}\overline{\lim\limits_{\varepsilon\to0}}\quad
\dashint_{B(z_0,\varepsilon)}K^T_{\mu}(z,z_0)\,dm(z)\ <\ \infty\
.\end{equation}

Then all the conclusions of Theorem 1 on solutions for the Dirichlet
problem (\ref{eqDIRICHLET}) with arbitrary continuous boundary data
$\varphi:\partial D\to{\Bbb R}$ to the Beltrami equation
(\ref{eqBeltrami}) with the source $\sigma$ hold.}

\medskip

Since $K^T_{\mu}(z,z_0) \leqslant K_{\mu}(z)$ for all $z$ and
$z_0\in \Bbb C$, we also obtain the following consequences of
Theorem 1 with the BMO type criterion.

\medskip

{\bf Corollary 4.} {\it Let $D$ be a bounded simply connected domain
in $\mathbb C$, $\sigma\in L_p(D)$, $p>2$, with compact support,
$\mu:D\to{\Bbb C}$ be a measurable function with $|\mu(z)|<1$ a.e.,
$K_{\mu}$ be locally bounded in $D$ and $K_{\mu}$ have a
do\-mi\-nant $Q\in $ {\rm BMO}$_{\rm loc}$ in a neighborhood of
$\partial D$. Then the conclusions of Theorem 1 hold.}

\medskip

{\bf Remark 5.} In particular, the conclusions of Theorem 1 hold if
$Q\in{\rm W}^{1,2}_{\rm loc}$ in a neighborhood of $\partial D$,
because of $W^{\,1,2}_{\rm loc} \subset {\rm VMO}_{\rm loc}$, see
e.g. \cite{BN}.

\medskip

{\bf Corollary 5.} {\it Let $D$ be a bounded simply connected domain
in $\mathbb C$, $\sigma\in L_p(D)$, $p>2$, with compact support,
$\mu:D\to{\Bbb C}$ be a measurable function with $|\mu(z)|<1$ a.e.,
$K_{\mu}$ be locally bounded in $D$ and $K_{\mu}$ have a
do\-mi\-nant $Q\in $ {\rm FMO} in a neighborhood of $\partial D$.
Then the conclusions of Theorem 1 hold. }

\medskip

Similarly, choosing in Lemma 1 the function $\psi(t)=1/t$, see also
Remark 1, we come to the next statement with the Calderon-Zygmund
type criterion.

\medskip

{\bf Theorem 2.} {\it Let $D$ be a bounded simply connected domain
in $\mathbb C$, $\sigma\in L_p(D)$, $p>2$, with compact support,
 $\mu:D\to{\Bbb C}$ be a measurable function with $|\mu(z)|<1$
a.e., $K_{\mu}$ be locally bounded in $D$, $K_{\mu}\in L_1(D)$ and,
for each point $z_0\in\partial{D}$ and
$\varepsilon_0=\varepsilon(z_0)>0$,
\begin{equation}\label{eqLOG}
\int\limits_{\varepsilon<|z-z_0|<\varepsilon_0}K^T_{\mu}(z,z_0)\,\frac{dm(z)}{|z-z_0|^2}
=o\left(\left[\log\frac{1}{\varepsilon}\right]^2\right)\ \ \ \hbox{
as $\varepsilon\to 0$}\ .\end{equation}

Then the Beltrami equation (\ref{eqBeltrami}) with the source
$\sigma$ has a locally H\"older continuous solution $\omega$ in the
class $W^{1,2}_{\rm loc}$ of the Dirichlet problem
(\ref{eqDIRICHLET}) in $D$ for each continuous function
$\varphi:\partial D\to{\Bbb R}$ that is unique up to an additive
pure imaginary constant.

Moreover, $\omega =h\circ f$, $h := {\cal A} + H$, where $f:\mathbb
C\to\mathbb C$ is a $\mu -$conformal mapping with $\mu$ extended by
zero outside of $D$, $H:D_*\to\mathbb C$  is a generalized analytic
function in $D_*:=f(D)$ with the source $S$ calculated in (\ref{r})
and ${\cal A}$ is a holomorphic function in $D_*$ with the Dirichlet
condition (\ref{eqSTAR}). }

\medskip

{\bf Remark 6.} Choosing in Lemma 1 the function
$\psi(t)=1/(t\log{1/t})$ instead of $\psi(t)=1/t$, we are able to
replace (\ref{eqLOG}) by the conditions
\begin{equation}\label{eqLOGLOG}
\int\limits_{\varepsilon<|z-z_0|<\varepsilon_0}\frac{K^T_{\mu}(z,z_0)\,dm(z)}
{\left(|z-z_0|\log{\frac{1}{|z-z_0|}}\right)^2}
=o\left(\left[\log\log\frac{1}{\varepsilon}\right]^2\right)\qquad\forall\
z_0\in\partial{D}\end{equation} as $\varepsilon\to 0$ for some
$\varepsilon_0=\varepsilon(z_0)>0$. More generally, we would be able
to give here the whole scale of the corresponding conditions in
$\log$ using functions $\psi(t)$ of the form
$1/(t\log{1}/{t}\cdot\log\log{1}/{t}\cdot\ldots\cdot\log\ldots\log{1}/{t})$.

\medskip

Choosing in Lemma 1 the functional parameter ${\psi}_{z_0}(t) : =
1/[tk^T_{\mu}(z_0,t)]$, where $k_{\mu}^T(z_0,r)$ is the integral
mean of $K^T_{{\mu}}(z,z_0)$ over the circle $S(z_0,r)\, :=\, \{ z
\in\mathbb C:\, |z-z_0|\,=\, r\}$, we obtain the Lehto type
criterion.

\medskip

{\bf Theorem 3.} {\it Let $D$ be a bounded simply connected domain
in $\mathbb C$, $\sigma\in L_p(D)$, $p>2$, with compact support,
$\mu:D\to{\Bbb C}$ be a measurable function with $|\mu(z)|<1$ a.e.,
$K_{\mu}$ be locally bounded in $D$, $K_{\mu}\in L_1(D)$ and, for
each point $z_0\in\partial{D}$ and
$\varepsilon_0=\varepsilon(z_0)>0$,
\begin{equation}\label{eqLEHTO}\int\limits_{0}^{\varepsilon_0}
\frac{dr}{rk^T_{\mu}(z_0,r)}\ =\ \infty\ .\end{equation}

Then the Beltrami equation (\ref{eqBeltrami}) with the source
$\sigma$ has a locally H\"older continuous solution $\omega$ in the
class $W^{1,2}_{\rm loc}$ of the Dirichlet problem
(\ref{eqDIRICHLET}) in $D$ for each continuous function
$\varphi:\partial D\to{\Bbb R}$ that is unique up to an additive
pure imaginary constant.

Moreover, $\omega =h\circ f$, $h := {\cal A} + H$, where $f:\mathbb
C\to\mathbb C$ is a $\mu -$conformal mapping with $\mu$ extended by
zero outside of $D$, $H:D_*\to\mathbb C$  is a generalized analytic
function in $D_*:=f(D)$ with the source $S$ calculated in (\ref{r})
and ${\cal A}$ is a holomorphic function in $D_*$ with the Dirichlet
condition (\ref{eqSTAR}). }

\medskip

{\bf Corollary 6.} {\it Let $D$ be a bounded simply connected domain
in $\mathbb C$, $\sigma\in L_p(D)$, $p>2$, with compact support,
$\mu:D\to{\Bbb C}$ be a measurable function with $|\mu(z)|<1$ a.e.,
$K_{\mu}$ be locally bounded in $D$, $K_{\mu}\in L_1(D)$ and, for
each point $z_0\in\partial{D}$,
\begin{equation}\label{eqLOGk}k^T_{\mu}(z_0,\varepsilon)=O\left(\log\frac{1}{\varepsilon}\right)
\qquad\mbox{as}\ \varepsilon\to0\ .\end{equation}

Then all conclusions of Theorem 3 on solutions for the Dirichlet
problem (\ref{eqDIRICHLET}) with arbitrary continuous boundary data
$\varphi:\partial D\to{\Bbb R}$ to the Beltrami equation
(\ref{eqBeltrami}) with the source $\sigma$ hold. }

\medskip

{\bf Remark 7.} In particular, the conclusions of Theorem 3 hold if
\begin{equation}\label{eqLOGK} K^T_{\mu}(z,z_0)=O\left(\log\frac{1}{|z-z_0|}\right)\qquad{\rm
as}\quad z\to z_0\quad\forall\ z_0\in\partial{D}\,.\end{equation}
Moreover, the condition (\ref{eqLOGk}) can be replaced by the series
of weaker conditions
\begin{equation}\label{edLOGLOGk}
k^T_{\mu}(z_0,\varepsilon)=O\left(\left[\log\frac{1}{\varepsilon}\cdot\log\log\frac{1}
{\varepsilon}\cdot\ldots\cdot\log\ldots\log\frac{1}{\varepsilon}
\right]\right) \qquad\forall\ z_0\in \partial{D}\ .
\end{equation}


To get another criterion, we need a couple of auxiliary statements.
The first of them can be found e.g. as Theorem 3.2 in
\cite{RSY$_5$}.

\medskip

{\bf Proposition 3.} {\it Let $Q:{\Bbb D}\to[0,\infty]$ be a
measurable function such that
\begin{equation}\label{eq5.555} \int\limits_{\Bbb
D}\Phi(Q(z))\,dm(z)<\infty\end{equation} where
$\Phi:[0,\infty]\to[0,\infty]$ is a non-decreasing convex function
such that \begin{equation}\label{eq3.333a}
\int\limits_{\delta}^{\infty}\frac{d\tau}{\tau\Phi^{-1}(\tau)}=\infty\end{equation}
for some $\delta>\Phi(+0)$. Then \begin{equation}\label{eq3.333A}
\int\limits_0^1\frac{dr}{rq(r)}=\infty\end{equation} where $q(r)$ is
the average of the function $Q(z)$ over the circle $|z|=r$.}

\medskip

Above we used the following notions of the inverse function for
monotone functions. Namely, for every non-decreasing function
$\Phi:[0,\infty]\to[0,\infty]$ the inverse function
$\Phi^{-1}:[0,\infty]\to[0,\infty]$ can be well-defined by setting
\begin{equation}\label{eqINVERSE}
\Phi^{-1}(\tau)\ :=\ \inf\limits_{\Phi(t)\geq\tau} t
\end{equation}
Here $\inf$ is equal to $\infty$ if the set of $t\in[0,\infty]$ such
that $\Phi(t)\geq\tau$ is empty. Note that the function $\Phi^{-1}$
is non-decreasing, too. It is also evident immediately by the
definition that $\Phi^{-1}(\Phi(t)) \leq t$ for all $t\in[0,\infty]$
with the equality except intervals of constancy of the function
$\Phi(t)$.

\medskip

Recall connections between integral conditions, see e.g. Theorem 2.5
in \cite{RSY$_5$}.

\medskip

{\bf Proposition 4.} Let $\Phi:[0,\infty]\to[0,\infty]$ be a
non-decreasing function and set
\begin{equation}\label{eqLOGFi}
H(t)\ =\ \log\Phi(t)\ .
\end{equation}
Then the equality
\begin{equation}\label{eq333Frer}\int\limits_{\Delta}^{\infty}H'(t)\,\frac{dt}{t}=\infty,
\end{equation}
implies the equality
\begin{equation}\label{eq333F}\int\limits_{\Delta}^{\infty}
\frac{dH(t)}{t}=\infty\,,\end{equation} and (\ref{eq333F}) is
equivalent to
\begin{equation}\label{eq333B}
\int\limits_{\Delta}^{\infty}H(t)\,\frac{dt}{t^2}=\infty\,\end{equation}
for some $\Delta>0$, and (\ref{eq333B}) is equivalent to each of the
equalities
\begin{equation}\label{eq333C} \int\limits_{0}^{\delta_*}H\left(\frac{1}{t}\right)\,{dt}=\infty\end{equation} for
some $\delta_*>0$, \begin{equation}\label{eq333D}
\int\limits_{\Delta_*}^{\infty}\frac{d\eta}{H^{-1}(\eta)}=\infty\end{equation}
for some $\Delta_*>H(+0)$ and to (\ref{eq3.333a}) for some
$\delta>\Phi(+0)$.

Moreover, (\ref{eq333Frer}) is equivalent to (\ref{eq333F}) and to
hence (\ref{eq333Frer})–(\ref{eq333D}) as well as to
(\ref{eq3.333a}) are equivalent to each other if $\Phi$ is in
addition absolutely continuous. In particular, all the given
conditions are equivalent if $\Phi$ is convex and non-decreasing.

\medskip

Note that the integral in (\ref{eq333F}) is understood as the
Lebesgue--Stieltjes integral and the integrals in (\ref{eq333Frer})
and (\ref{eq333B})--(\ref{eq333D}) as the ordinary Lebesgue
integrals. It is necessary to give one more explanation. From the
right hand sides in the conditions (\ref{eq333Frer})--(\ref{eq333D})
we have in mind $+\infty$. If $\Phi(t)=0$ for $t\in[0,t_*$, then
$H(t)=-\infty$ for $t\in[0,t_*]$ and we complete the definition
$H'(t)=0$ for $t\in[0,t_*]$. Note, the conditions (\ref{eq333F}) and
(\ref{eq333B}) exclude that $t_*$ belongs to the interval of
integrability because in the contrary case the left hand sides in
(\ref{eq333F}) and (\ref{eq333B}) are either equal to $-\infty$ or
indeterminate. Hence we may assume in
(\ref{eq333Frer})--(\ref{eq333C}) that $\delta >t_0$,
correspondingly, $\Delta<1/t_0$ where
$t_0:=\sup\limits_{\Phi(t)=0}t$, and set $t_0=0$ if $\Phi(0)>0$. The
most interesting condition (\ref{eq333B}) can be written in the
form:
\begin{equation}\label{eq5!}
\int\limits_{\Delta}^{\infty}\log\, \Phi(t)\ \frac{dt}{t^{2}}\ =\
+\infty\ \ \ \ \ \ \mbox{for some $\Delta > 0$}\ .
\end{equation}

Combining Proposition 3 and 4 with Theorems 3 we obtain the
following significant result with the Orlicz type criterion.

\medskip

{\bf Theorem 4.} {\it Let $D$ be a bounded simply connected domain
in $\mathbb C$, $\sigma\in L_p(D)$, $p>2$, with compact support,
$\mu:D\to{\Bbb C}$ be a measurable function with $|\mu(z)|<1$ a.e.,
$K_{\mu}$ be locally bounded in $D$, $K_{\mu}\in L_1(D)$ and, for
each point $z_0\in
\partial{D}$ and a neighborhood $U_{z_0}$ of $z_0$,
\begin{equation}\label{eqINTEG}\int\limits_{U_{z_0}}\Phi_{z_0}\left(K^T_{\mu}(z,z_0)\right)\,dm(z)<\infty\
,
\end{equation} where $\Phi_{z_0}:(0,\infty]\to(0,\infty]$ is a convex non-decreasing function
such that
\begin{equation}\label{eqINT}
\int\limits_{\Delta(z_0)}^{\infty}\log\,\Phi_{z_0}(t)\,\frac{dt}{t^2}\
=\ +\infty\ \ \ \hbox{for some $\Delta(z_0)>0$}\ .\end{equation}

Then the Beltrami equation (\ref{eqBeltrami}) with the source
$\sigma$ has a locally H\"older continuous solution $\omega$ in the
class $W^{1,2}_{\rm loc}$ of the Dirichlet problem
(\ref{eqDIRICHLET}) in $D$ for each continuous function
$\varphi:\partial D\to{\Bbb R}$ that is unique up to an additive
pure imaginary constant.

Moreover, $\omega =h\circ f$, $h := {\cal A} + H$, where $f:\mathbb
C\to\mathbb C$ is a $\mu -$conformal mapping with $\mu$ extended by
zero outside of $D$, $H:D_*\to\mathbb C$  is a generalized analytic
function in $D_*:=f(D)$ with the source $S$ calculated in (\ref{r})
and ${\cal A}$ is a holomorphic function in $D_*$ with the Dirichlet
condition (\ref{eqSTAR}). }

\medskip

{\bf Corollary 7.} {\it Let $D$ be a bounded simply connected domain
in $\mathbb C$, $\sigma\in L_p(D)$, $p>2$, with compact support,
$\mu:D\to{\Bbb C}$ be a measurable function with $|\mu(z)|<1$ a.e.,
$K_{\mu}$ be locally bounded in $D$, $K_{\mu}\in L_1(D)$ and, for
each point $z_0\in \partial{D}$, a neighborhood $U_{z_0}$ of $z_0$
and $\alpha(z_0)>0$,
\begin{equation}\label{eqEXP}\int\limits_{U_{z_0}}e^{\alpha(z_0)
K^T_{\mu}(z,z_0)}\,dm(z)<\infty\ .
\end{equation}

Then all conclusions of Theorem 4 on solutions for the Dirichlet
problem (\ref{eqDIRICHLET}) with continuous data $\varphi:\partial
D\to{\Bbb R}$ to the Beltrami equation (\ref{eqBeltrami}) with the
source $\sigma$ hold. }

\medskip

{\bf Corollary 8.} {\it Let $D$ be a bounded simply connected domain
in $\mathbb C$, $\sigma\in L_p(D)$, $p>2$, with compact support,
$\mu:D\to{\Bbb C}$ be a measurable function with $|\mu(z)|<1$ a.e.,
$K_{\mu}$ be locally bounded in $D$ and, for a neighborhood $U$ of
$\partial D$,
\begin{equation}\label{eqINTK}\int\limits_{U}\Phi\left(K_{\mu}(z)\right)\,dm(z)<\infty\ ,\end{equation}
where $\Phi:(0,\infty]\to(0,\infty]$ is a convex non-decreasing
function with, for $\delta>0$,
\begin{equation}\label{eqINTF}
\int\limits_{\delta}^{\infty}\log\,\Phi(t)\,\frac{dt}{t^2}\ =\
+\infty\ .\end{equation}

Then all conclusions of Theorem 4 on solutions for the Dirichlet
problem (\ref{eqDIRICHLET}) with continuous data $\varphi:\partial
D\to{\Bbb R}$ to the Beltrami equation (\ref{eqBeltrami}) with the
source $\sigma$ hold. }

\medskip

{\bf Remark 8.} By Theorems 2.5 and 5.1 in \cite{RSY$_5$}, condition
(\ref{eqINTF}) is not only sufficient but also necessary to have the
regular solutions of the Dirichlet problem (\ref{eqDIRICHLET}) in
$D$ for arbitrary Beltrami equations with sources
(\ref{eqBeltrami}), satisfying the integral constraints
(\ref{eqINTK}), for all continuous functions $\varphi:\partial
D\to\Bbb{R}$ because such solutions have the representation through
regular homeomorphic solutions $f=f^{\mu}$ of the homogeneous
Beltrami equation (\ref{a}) from Proposition 1.

\medskip

{\bf Corollary 9.} {\it Let $D$ be a bounded simply connected domain
in $\mathbb C$, $\sigma\in L_p(D)$, $p>2$, with compact support,
$\mu:D\to{\Bbb C}$ be a measurable function with $|\mu(z)|<1$ a.e.,
$K_{\mu}$ be locally bounded in $D$ and, for a neighborhood $U$ of
$\partial D$ and $\alpha>0$,
\begin{equation}\label{eqEXPA}\int\limits_{U}e^{\alpha K_{\mu}(z)}\,dm(z)\ <\
\infty\ .
\end{equation}

Then all conclusions of Theorem 4 on solutions for the Dirichlet
problem (\ref{eqDIRICHLET}) with continuous data $\varphi:\partial
D\to{\Bbb R}$ to the Beltrami equation (\ref{eqBeltrami}) with the
source $\sigma$ hold. }

\section{The main lemma in general domains}

In this section we obtain criteria for the existence, representation
and regularity of the so-called multi-valued solutions $\omega$ of
the Dirichlet problem (\ref{eqDIRICHLET}) to the Beltrami equations
with sources (\ref{eqBeltrami}) in the spirit of the theory of
multi-valued analytic functions in arbitrary bounded domains $D$ in
$\mathbb C$ with no boundary component degenerated to a single
point. Simple examples show that such domains form the most wide
class of domains for which the problem is always solvable for any
continuous boundary data.

We say that a locally H\"older continuous function $\omega
:B(z_0,\varepsilon_0)\to{\Bbb C}$, where
$B(z_0,\varepsilon_0)\subseteq D$, is a {\bf local regular solution
of the equation} (\ref{eqBeltrami}) if $\omega\in W_{\rm loc}^{1,2}$
and $\omega$ satisfies (\ref{eqBeltrami}) a.e. in
$B(z_0,\varepsilon_0)$. Local regular solutions
$\omega_0:B(z_0,\varepsilon_0)\to{\Bbb C}$ and
$\omega_*:B(z_*,\varepsilon_*)\to{\Bbb C}$ of the equation
(\ref{eqBeltrami}) will be called extension of each to other if
there is a finite chain of its local regular solutions
$\omega_i:B(z_i,\varepsilon_i)\to\Bbb{C}$, $i=1,\ldots,m$, such that
$\omega_1=\omega_0$, $\omega_m=\omega_*$ and $\omega_i(z)=
\omega_{i+1}(z)$ for all points $z\in E_i:=B(z_i,\varepsilon_i)\cap
B(z_{i+1},\varepsilon_{i+1})\neq\emptyset$, $i=1,\ldots,m-1$.
\par
A collection of local regular solutions
$\omega_j:B(z_j,\varepsilon_j)\to{\Bbb C}$, $j\in J$, will be called
a {\bf regular multi-valued solution} of the equation
(\ref{eqBeltrami}) in $D$ if the collection of the disks
$B(z_j,\varepsilon_j)$ cover the domain $D$ and $\omega_j$ are
extensions of each to other through this collection and the
collection is maximal by inclusion.
\par
A regular multi-valued solution of the equation (\ref{eqBeltrami})
will be called a {\bf regular multi-valued solution of the Dirichlet
problem} (\ref{eqDIRICHLET}) to (\ref{eqBeltrami}) in $D$ if
$u(z)={\rm Re}\,\omega(z)={\rm Re}\,\omega_{j}(z)$, $z\in
B(z_j,\varepsilon_j)$, $j\in J$, is a single-valued function in $D$
satisfying the Dirichlet condition
$\lim\limits_{z\in\zeta}u(z)=\varphi(\zeta)$ for all $\zeta\in
\partial D$.




{\bf Lemma 3.} {\it Let $D$ be a bounded domain in ${\Bbb C}$ with
no boundary component degenerated to a single point, $\sigma \in
L_p(D)$, $p>2$, with compact support in $D$, $\mu:D\to{\Bbb C}$ be a
measurable function with $|\mu(z)|<1$ a.e., $K_{\mu}$ be locally
bounded in $D$, $K_{\mu}\in L_1(D)$ and, for each point
$z_0\in\partial{D}$ and $\varepsilon_0=\varepsilon(z_0)>0$,
\begin{equation}\label{M3omalM}
\int\limits_{\varepsilon<|z-z_0|<\varepsilon_0}
K^T_{\mu}(z,z_0)\cdot\psi^2_{z_0,\varepsilon}(|z-z_0|)\,dm(z)=o(I_{z_0}^{2}(\varepsilon))\quad\hbox{as
$\varepsilon\to0$}\ ,
\end{equation} where
$\psi_{z_0,\varepsilon}: (0,\varepsilon_0)\to(0,\infty)$ is a family
of measurable functions such that
\begin{equation}\label{Meq3.5.3M}
I_{z_0}(\varepsilon)\colon
=\int\limits_{\varepsilon}^{\varepsilon_0}
\psi_{z_0,\varepsilon}(t)\,dt<\infty\qquad\forall\
\varepsilon\in(0,\varepsilon_0)\,.\end{equation}

Then the Beltrami equation (\ref{eqBeltrami}) with the source
$\sigma$ has a regular multi-valued solution $\omega$ of the
Dirichlet problem (\ref{eqDIRICHLET}) in $D$ for each continuous
function $\varphi:\partial D\to{\Bbb R}$ that is unique up to an
additive pure imaginary constant.

Moreover, $\omega =h\circ f$, $h := {\cal A} + H$, where $f:\mathbb
C\to\mathbb C$ is a $\mu -$conformal mapping with $\mu$ extended by
zero outside of $D$, $H:D_*\to\mathbb C$  is a generalized analytic
function in $D_*:=f(D)$ with the source $S$ calculated in (\ref{r})
and ${\cal A}$ is a multi-valued analytic function in $D_*$ with a
single valued real part satisfying the Dirichlet condition
(\ref{eqSTAR}). }

\begin{proof} By Proposition 1 there is a $\mu -$conformal
mapping $f:\mathbb C\to\mathbb C$ with $\mu$ extended by zero
outside of $D$, which is locally quasiconformal in $D$. Arguing
locally as in the first item of the proof to Lemma 1, we first show
that the desired solution is unique up to an additive pure imaginary
constant. Moreover, by the second item of the proof to Lemma 1, the
function $S$ described in (\ref{r}) belongs to the class
$L_{p_*}(D_*)$ in the domain $D_*=f(D)$ with $p_*\in(2,p)$.

Note also that $D_*$ is also a bounded domain with no boundary
components degenerated to a single point. Let
$\varphi^*:=\left(\varphi\circ f^{-1} - {\rm Re}\,
H\right)|_{\partial D_*}$, where $H$ is the generalized analytic
function (\ref{eqINTEGRAL}) with the source $S$. Then by Corollary
4.1.8 and Theorem 4.2.2 in \cite{Rans} there is a harmonic function
$u:D_*\to\mathbb R$ satisfying the Dirichlet condition
\begin{equation}\label{eqHARMONIC}
\lim_{w\to\zeta}\ u(w)\ =\ \varphi^*(\zeta)\ \ \ \ \ \forall\
\zeta\in\partial D_*\ .
\end{equation}

Now, let $B_0=B(z_0,r_0)$ be a disk in the domain $D$. Then ${\frak
B}_0 := f(B_0)$ is a simply connected subdomain of the domain $D_*$
where there is a conjugate harmonic function $v$ determined up to an
additive constant such that $u+iv$ is a single--valued analytic
function. Let us denote through ${\cal A}_0$ the holomorphic
function corresponding to the choice of such a harmonic function
$v_0$ in ${\frak B}_0$ with the normalization $v_0(f(z_0))=0$.
Thereby we have determined the initial element of a multi-valued
analytic function. The function ${\cal A}_0$ can be extended to,
generally speaking multi-valued, analytic function ${\cal A}$  along
any path in $D_*$ because $u$ is given in the whole domain  $D_*$.

Thus, $\omega :=h\circ f$, $h={\cal A}+H$, is a continuous
multi-valued solution of the Dirichlet problem (\ref{eqDIRICHLET})
in $D$ for the Beltrami equation with the source $\sigma$
(\ref{eqBeltrami}) of the class $W^{1,1}_{\rm loc}$ by Lemma III.6.4
in \cite{LV}, see also Remark 1. Arguing as in the second item of
the proof to Lemma 1 with the application of the auxiliary
quasiconformal mapping $f^{\tilde\mu}:\mathbb C\to\mathbb C$, it is
easy to prove on the basis of Lemma 1 in \cite{GNRY*} that
$\omega\in W^{1,2}_{\rm loc}$. Finally, the given solution $\omega$
is locally H\"older continuous because the function $h$ and the
mapping $f$ are so, see Remark 1.
\end{proof}

{\bf Remark 9.} Note that if the family of the functions
$\psi_{z_0,\varepsilon}(t)\equiv\psi_{z_0}(t)$, $z_0\in\partial{D}$,
in Lemma 3 is independent on the parameter $\varepsilon$, then the
condition (\ref{M3omalM}) implies that $I_{z_0}(\varepsilon)\to
\infty$ as $\varepsilon\to 0$. This follows immediately from
arguments by contradiction, apply for it (\ref{eqConnect}) and the
condition $K_{\mu}\in L_1(D)$. Note also that (\ref{M3omalM}) holds,
in particular, if, for some $\varepsilon_0=\varepsilon(z_0)$,
\begin{equation}\label{M333omal}
\int\limits_{|z-z_0|<\varepsilon_0}
K^T_{\mu}(z,z_0)\cdot\psi_{z_0}^2(|z-z_0|)\,dm(z)<\infty \qquad
\forall\ z_0\in\partial{D}\end{equation} and
$I_{z_0}(\varepsilon)\to \infty$ as $\varepsilon\to 0$. In other
words, for the existence of  a regular multi-valued solutions of the
Dirichlet problem (\ref{eqDIRICHLET}) in $D$ for the Beltrami
equations with sources (\ref{eqBeltrami}) for all con\-ti\-nu\-ous
boundary functions $\varphi$, it is sufficient that the integral in
(\ref{M333omal}) converges for some nonnegative function
$\psi_{z_0}(t)$ that is locally integrable over $(0,\varepsilon_0 ]$
but has a nonintegrable singularity at $0$. The functions
$\log^{\lambda}(e/|z-z_0|)$, $\lambda\in (0,1)$, $z\in\Bbb D$,
$z_0\in\overline{\Bbb D}$, and $\psi(t)=1/(t \,\, \log(e/t))$,
$t\in(0,1)$, show that the condition (\ref{M333omal}) is compatible
with the condition $I_{z_0}(\varepsilon)\to\infty $ as
$\varepsilon\to 0$. Furthermore, the condition (\ref{M3omalM}) shows
that it is sufficient for the existence of  regular multi-valued
solutions of the Dirichlet problem (\ref{eqDIRICHLET}) in $D$ for
the Beltrami equations with sources (\ref{eqBeltrami}) for all
con\-ti\-nu\-ous boundary functions $\varphi$ even if the integral
in (\ref{M333omal}) is divergent in a controlled way.

\section{The main criteria in general domains}

Arguing as in Section 3, we derive from Lemma 3 the following
consequences.

\medskip

{\bf Theorem 5.} {\it Let $D$ be a bounded domain in ${\Bbb C}$ with
no boundary component degenerated to a single point, $\sigma \in
L_p(D)$, $p>2$, have compact support in $D$, $\mu:D\to{\Bbb C}$ be a
measurable function with $|\mu(z)|<1$ a.e., $K_{\mu}$ be locally
bounded in $D$, $K_{\mu}\in L_1(D)$, $K^T_{\mu}(z,z_0)\leqslant
Q_{z_0}(z)$ a.e. in $U_{z_0}$ for each point $z_0\in \partial{D}$, a
neighborhood $U_{z_0}$ of $z_0$ and a function $Q_{z_0}:
U_{z_0}\to[0,\infty]$ in the class ${\rm FMO}({z_0})$.

Then the Beltrami equation (\ref{eqBeltrami}) with the source
$\sigma$ has a regular multi-valued solution $\omega$ of the
Dirichlet problem (\ref{eqDIRICHLET}) in $D$ for each continuous
function $\varphi:\partial D\to{\Bbb R}$ that is unique up to an
additive pure imaginary constant.

Moreover, $\omega =h\circ f$, $h := {\cal A} + H$, where $f:\mathbb
C\to\mathbb C$ is a $\mu -$conformal mapping with $\mu$ extended by
zero outside of $D$, $H:D_*\to\mathbb C$  is a generalized analytic
function in $D_*:=f(D)$ with the source $S$ calculated in (\ref{r})
and ${\cal A}$ is a multi-valued analytic function in $D_*$ with a
single valued real part satisfying the Dirichlet condition
(\ref{eqSTAR}).}

\medskip

{\bf Corollary 10.} {\it Let $D$ be a bounded domain in ${\Bbb C}$
with no boundary component degenerated to a single point, $\sigma
\in L_p(D)$, $p>2$, have compact support in $D$,
 $\mu:D\to{\Bbb C}$ be a measurable function with $|\mu(z)|<1$
a.e., $K_{\mu}$ be locally bounded in $D$, $K_{\mu}\in L_1(D)$ and,
for each point $z_0\in\partial{D}$,
\begin{equation}\label{MeqMEAN}\overline{\lim\limits_{\varepsilon\to0}}\quad
\dashint_{B(z_0,\varepsilon)}K^T_{\mu}(z,z_0)\,dm(z)<\infty\
.\end{equation} Then all the conclusions of Theorem 5 hold.}

\medskip

{\bf Corollary 11.} {\it Let $D$ be a bounded domain in ${\Bbb C}$
with no boundary component degenerated to a single point, $\sigma
\in L_p(D)$, $p>2$, have compact support in $D$, $\mu:D\to{\Bbb C}$
be a measurable function with $|\mu(z)|<1$ a.e., $K_{\mu}$ be
locally bounded in $D$ and $K_{\mu}$ have a do\-mi\-nant $Q\in $
{\rm BMO}$_{\rm loc}$ in a neighborhood of $\partial D$. Then all
the conclusions of Theorem 5 hold.}

\medskip

{\bf Remark 10.} In particular, the conclusions of Theorem 5 hold if
$Q\in{\rm W}^{1,2}_{\rm loc}$ in a neighborhood of $\partial D$,
because of $W^{\,1,2}_{\rm loc} \subset {\rm VMO}_{\rm loc}$.

\medskip

{\bf Corollary 12.} {\it Let $D$ be a bounded domain in ${\Bbb C}$
with no boundary component degenerated to a single point, $\sigma
\in L_p(D)$, $p>2$, have compact support in $D$, $\mu:D\to{\Bbb C}$
be a measurable function with $|\mu(z)|<1$ a.e., $K_{\mu}$ be
locally bounded in $D$ and $K_{\mu}$ have a do\-mi\-nant $Q\in $
{\rm FMO} in a neighborhood of $\partial D$. Then all the
conclusions of Theorem 5 hold. }

\medskip

{\bf Theorem 6.} {\it Let $D$ be a bounded domain in ${\Bbb C}$ with
no boundary component degenerated to a single point, $\sigma \in
L_p(D)$, $p>2$, have compact support in $D$, $\mu:D\to{\Bbb C}$ be a
measurable function with $|\mu(z)|<1$ a.e., $K_{\mu}$ be locally
bounded in $D$, $K_{\mu}\in L_1(D)$ and, for each point
$z_0\in\partial{D}$ and $\varepsilon_0=\varepsilon(z_0)>0$,
\begin{equation}\label{MeqLOG}
\int\limits_{\varepsilon<|z-z_0|<\varepsilon_0}K^T_{\mu}(z,z_0)\,\frac{dm(z)}{|z-z_0|^2}
=o\left(\left[\log\frac{1}{\varepsilon}\right]^2\right)\qquad\hbox{as
$\varepsilon\to 0$}\ .\end{equation}

Then the Beltrami equation (\ref{eqBeltrami}) with the source
$\sigma$ has a regular multi-valued solution $\omega$ of the
Dirichlet problem (\ref{eqDIRICHLET}) in $D$ for each continuous
function $\varphi:\partial D\to{\Bbb R}$ that is unique up to an
additive pure imaginary constant.

Moreover, $\omega =h\circ f$, $h := {\cal A} + H$, where $f:\mathbb
C\to\mathbb C$ is a $\mu -$conformal mapping with $\mu$ extended by
zero outside of $D$, $H:D_*\to\mathbb C$  is a generalized analytic
function in $D_*:=f(D)$ with the source $S$ calculated in (\ref{r})
and ${\cal A}$ is  a multi-valued analytic function in $D_*$ with a
single valued real part satisfying the Dirichlet condition
(\ref{eqSTAR}). }

\medskip

{\bf Remark 11.} Choosing in Lemma 3 the function
$\psi(t)=1/(t\log{1/t})$ instead of $\psi(t)=1/t$, we are able to
replace (\ref{MeqLOG}) by the conditions
\begin{equation}\label{MeqLOGLOG}
\int\limits_{\varepsilon<|z-z_0|<\varepsilon_0}\frac{K^T_{\mu}(z,z_0)\,dm(z)}
{\left(|z-z_0|\log{\frac{1}{|z-z_0|}}\right)^2}
=o\left(\left[\log\log\frac{1}{\varepsilon}\right]^2\right)\qquad\forall\
z_0\in\partial{D}\end{equation} as $\varepsilon\to 0$ for some
$\varepsilon_0=\varepsilon(z_0)>0$. More generally, we would be able
to give here the whole scale of the corresponding conditions in
$\log$ using functions $\psi(t)$ of the form
$1/(t\log{1}/{t}\cdot\log\log{1}/{t}\cdot\ldots\cdot\log\ldots\log{1}/{t})$.

\medskip

{\bf Theorem 7.} {\it Let $D$ be a bounded domain in ${\Bbb C}$ with
no boundary component degenerated to a single point, $\sigma \in
L_p(D)$, $p>2$, have compact support in $D$, $\mu:D\to{\Bbb C}$ be a
measurable function with $|\mu(z)|<1$ a.e., $K_{\mu}$ be locally
bounded in $D$, $K_{\mu}\in L_1(D)$ and, for each point
$z_0\in\partial{D}$ and for some $\varepsilon_0=\varepsilon(z_0)>0$,
\begin{equation}\label{MeqLEHTO}\int\limits_{0}^{\varepsilon_0}
\frac{dr}{rk^T_{\mu}(z_0,r)}\ =\ \infty\ ,\end{equation} where
$k^T_{\mu}(z_0,r)$ is the mean value of $K^T_{\mu}(z_0,r)$ over the
circles $S(z_0, r)$.

Then the Beltrami equation (\ref{eqBeltrami}) with the source
$\sigma$ has a regular multi-valued solution $\omega$ of the
Dirichlet problem (\ref{eqDIRICHLET}) in $D$ for each continuous
function $\varphi:\partial D\to{\Bbb R}$ that is unique up to an
additive pure imaginary constant.

Moreover, $\omega =h\circ f$, $h := {\cal A} + H$, where $f:\mathbb
C\to\mathbb C$ is a $\mu -$conformal mapping with $\mu$ extended by
zero outside of $D$, $H:D_*\to\mathbb C$  is a generalized analytic
function in $D_*:=f(D)$ with the source $S$ calculated in (\ref{r})
and ${\cal A}$ is  a multi-valued analytic function in $D_*$ with a
single valued real part satisfying the Dirichlet condition
(\ref{eqSTAR}). }

\medskip

{\bf Corollary 13.} {\it Let $D$ be a bounded domain in ${\Bbb C}$
with no boundary component degenerated to a single point, $\sigma
\in L_p(D)$, $p>2$, have compact support in $D$, $\mu:D\to{\Bbb C}$
be a measurable function with $|\mu(z)|<1$ a.e., $K_{\mu}$ be
locally bounded in $D$, $K_{\mu}\in L_1(D)$ and, for each point
$z_0\in\partial{D}$,
\begin{equation}\label{MeqLOGk}k^T_{\mu}(z_0,\varepsilon)=O\left(\log\frac{1}{\varepsilon}\right)
\qquad\mbox{as}\ \varepsilon\to0\ .\end{equation}

Then all conclusions of Theorem 7 on regular multi-valued solutions
for the Dirichlet problem (\ref{eqDIRICHLET}) with arbitrary
continuous boundary data $\varphi:\partial D\to{\Bbb R}$ to the
Beltrami equation (\ref{eqBeltrami}) with the source $\sigma$ hold.
}

\medskip

{\bf Remark 12.} In particular, the conclusions of Theorem 7 hold if
\begin{equation}\label{eqLOGK} K^T_{\mu}(z,z_0)=O\left(\log\frac{1}{|z-z_0|}\right)\qquad{\rm
as}\quad z\to z_0\quad\forall\ z_0\in\partial{D}\,.\end{equation}
Moreover, the condition (\ref{MeqLOGk}) can be replaced by the
series of weaker conditions
\begin{equation}\label{MedLOGLOGk}
k^T_{\mu}(z_0,\varepsilon)=O\left(\left[\log\frac{1}{\varepsilon}\cdot\log\log\frac{1}
{\varepsilon}\cdot\ldots\cdot\log\ldots\log\frac{1}{\varepsilon}
\right]\right) \qquad\forall\ z_0\in \partial{D}\ .
\end{equation}

\bigskip


{\bf Theorem 8.} {\it Let $D$ be a bounded domain in ${\Bbb C}$ with
no boundary component degenerated to a single point, $\sigma \in
L_p(D)$, $p>2$, have compact support in $D$, $\mu:D\to{\Bbb C}$ be a
measurable function with $|\mu(z)|<1$ a.e., $K_{\mu}$ be locally
bounded in $D$, $K_{\mu}\in L_1(D)$ and, for each point
$z_0\in\partial{D}$ and a neighborhood $U_{z_0}$ of $z_0$,
\begin{equation}\label{eqINTEG}\int\limits_{U_{z_0}}\Phi_{z_0}\left(K^T_{\mu}(z,z_0)\right)\,dm(z)<\infty
\ ,\end{equation} where $\Phi_{z_0}:(0,\infty]\to(0,\infty]$ is a
convex non-decreasing function such that
\begin{equation}\label{MeqINT}
\int\limits_{\Delta(z_0)}^{\infty}\log\,\Phi_{z_0}(t)\,\frac{dt}{t^2}\
=\ +\infty\ \ \ \hbox{for some $\Delta(z_0)>0$}\ .\end{equation}

Then the Beltrami equation (\ref{eqBeltrami}) with the source
$\sigma$ has a regular multi-valued solution $\omega$ of the
Dirichlet problem (\ref{eqDIRICHLET}) in $D$ for each continuous
function $\varphi:\partial D\to{\Bbb R}$ that is unique up to an
additive pure imaginary constant.

Moreover, $\omega =h\circ f$, $h := {\cal A} + H$, where $f:\mathbb
C\to\mathbb C$ is a $\mu -$conformal mapping with $\mu$ extended by
zero outside of $D$, $H:D_*\to\mathbb C$  is a generalized analytic
function in $D_*:=f(D)$ with the source $S$ calculated in (\ref{r})
and ${\cal A}$ is a multi-valued analytic function in $D_*$ with a
single valued real part satisfying the Dirichlet condition
(\ref{eqSTAR}). }

\bigskip

{\bf Corollary 14.} {\it Let $D$ be a bounded domain in ${\Bbb C}$
with no boundary component degenerated to a single point, $\sigma
\in L_p(D)$, $p>2$, have compact support in $D$, $\mu:D\to{\Bbb C}$
be a measurable function with $|\mu(z)|<1$ a.e., $K_{\mu}$ be
locally bounded in $D$, $K_{\mu}\in L_1(D)$ and, for each point
$z_0\in\partial{D}$ and a neighborhood $U_{z_0}$ of $z_0$,
\begin{equation}\label{MeqEXP}\int\limits_{U_{z_0}}e^{\alpha(z_0) K^T_{\mu}(z,z_0)}\,dm(z)<\infty
\qquad\hbox{for some $\alpha(z_0)>0$}\ .\end{equation}

Then all the conclusions of Theorem 8 on regular multi-valued
solutions for the Dirichlet problem (\ref{eqDIRICHLET}) with
continuous data $\varphi:\partial D\to{\Bbb R}$ to the Beltrami
equation (\ref{eqBeltrami}) with the source $\sigma$ hold. }

\bigskip


{\bf Corollary 15.} {\it Let $D$ be a bounded domain in ${\Bbb C}$
with no boundary component degenerated to a single point, $\sigma
\in L_p(D)$, $p>2$, have compact support in $D$, $\mu:D\to{\Bbb C}$
be a measurable function with $|\mu(z)|<1$ a.e., $K_{\mu}$ be
locally bounded in $D$ and, for a neighborhood $U$ of $\partial D$,
\begin{equation}\label{MeqINTK}\int\limits_{U}\Phi\left(K_{\mu}(z)\right)\,dm(z)<\infty\ ,\end{equation}
where $\Phi:(0,\infty]\to(0,\infty]$ is a convex non-decreasing
function with, for $\delta>0$,
\begin{equation}\label{MeqINTF}
\int\limits_{\delta}^{\infty}\log\,\Phi(t)\,\frac{dt}{t^2}\ =\
+\infty\ .\end{equation}

Then all the conclusions of Theorem 8 on regular multi-valued
solutions for the Dirichlet problem (\ref{eqDIRICHLET}) with
continuous data $\varphi:\partial D\to{\Bbb R}$ to the Beltrami
equation (\ref{eqBeltrami}) with the source $\sigma$ hold. }

\medskip

{\bf Remark 13.} The condition (\ref{MeqINTF}) is not only
sufficient but also ne\-ces\-sa\-ry to have the regular multi-valued
solutions of the Dirichlet problem (\ref{eqDIRICHLET}) for arbitrary
Beltrami equations with sources (\ref{eqBeltrami}), satisfying the
integral constraints (\ref{MeqINTK}), for all continuous functions
$\varphi:\partial D\to\Bbb{R}$, see arguments in Remark 8.

\medskip

{\bf Corollary 16.} {\it Let $D$ be a bounded domain in ${\Bbb C}$
with no boundary component degenerated to a single point, $\sigma
\in L_p(D)$, $p>2$, have compact support in $D$, $\mu:D\to{\Bbb C}$
be a measurable function with $|\mu(z)|<1$ a.e., $K_{\mu}$ be
locally bounded in $D$ and, for a neighborhood $U$ of $\partial D$
and some $\alpha>0$,
\begin{equation}\label{eqEXPA}\int\limits_{U}e^{\alpha K_{\mu}(z)}\,dm(z)\ <\
\infty\ .
\end{equation}

Then all the conclusions of Theorem 8 on regular multivalued
solutions for the Dirichlet problem (\ref{eqDIRICHLET}) with
continuous data $\varphi:\partial D\to{\Bbb R}$ to the Beltrami
equation (\ref{eqBeltrami}) with the source $\sigma$ hold. }

\section{Dirichlet problem for Poisson type equations}

Let us denote by $\mathbb S^{2\times 2}$ the collection of all
$2\times 2$ matrices with real entries
\begin{equation}\label{matrix}
A\  =\ \left[\begin{array}{ccc} a_{11}  & a_{12} \\
a_{21} & a_{22} \end{array}\right]\ ,\end{equation} which are
symmetric, i.e., $a_{12}=a_{21}$, with ${\rm det}\,A=1$ and {\bf
ellipticity condition} ${\rm det}\,(I+A)>0$, where $I$ is the unit
$2\times 2$ matrix. The latter condition means in terms of entries
of $A$ that $(1+a_{11})(1+a_{22})>a_{12}a_{21}$.


Now, let us consider in a domain $D$ of the complex plane $\mathbb
C$ the {Poisson type equations} (\ref{eqPotential}), where $A: D\to
\mathbb S^{2\times 2}$ is a measurable matrix valued function whose
elements $a_{ij}(z)$, $i,j=1,2$ are measurable and locally bounded
and first the source $g:D\to\mathbb R$ is a scalar function in
$L_{1,\rm loc}$.

\medskip

It is well--known, see Theorem 16.1.6 in \cite{AIM}, that
nonhomogeneous Beltrami equations (\ref{eqBeltrami}) with locally
bounded $K_{\mu}$ are closely connected with the Poisson type
equations (\ref{eqPotential}), where $A:D\to\mathbb S^{2\times 2}$
is the measurable matrix valued function
\begin{equation}\label{6.9}
A(z)\ :=\ \left[\begin{array}{ccc} {|1-\mu(z)|^2\over 1-|\mu(z)|^2}  & {-2{\rm Im}\,\mu(z)\over 1-|\mu(z)|^2} \\
                            {-2{\rm Im}\,\mu(z)\over 1-|\mu(z)|^2}          & {|1+\mu(z)|^2\over 1-|\mu(z)|^2}
                            \end{array}\right]\ ,
\end{equation}
whose entries $a_{ij}(z)$ are dominated by $K_{\mu}(z)$ and, thus,
they are locally bounded.


Vice versa, locally uniform elliptic (\ref{eqPotential}) with
measurable $A:D\to\mathbb S^{2\times 2}$ just correspond to
nonhomogeneous Beltrami equations (\ref{eqBeltrami}) with
coefficients
\begin{equation}\label{6.10}
\mu_A\ :=\ -\frac{a_{11}-a_{22}+i(a_{12}+a_{21})}{2+a_{11}+a_{22}}\
,
\end{equation} whose dilatation quotients $K_{\mu_A}$ are locally bounded.

Given such a matrix function $A$ and a $\mu $-conformal mapping
$f^{\mu }:D\rightarrow \mathbb{C}$, we have already seen in Lemma 1
of \cite{GNR-}, by direct computations, that if a function $T$ and
the entries of $A$ are sufficiently smooth, then
\begin{equation}
\mbox{div}\,[A(z)\nabla \,(T(f^{\mu }(z)))]\ =\ J(z)\triangle
\,T(f^{\mu }(z))\ .  \label{6.10k}
\end{equation}%
In the case $T\in W_{\mathrm{loc}}^{1,2}$, we understand equality (\ref%
{6.10k}) in the distributional sense, see Proposition 3.1 in
\cite{GNR}, i.e., for all $\psi \in W_{0}^{1,2}(D)$,
\begin{equation}
\int\limits_{D}\langle A\nabla (T\circ f^{\mu }),\nabla \psi \rangle
\ dm_{z}\ =\ \int\limits_{D}J(z)\langle M^{-1}((\nabla T)\circ
f^{\mu }),\nabla \psi \rangle \ dm_{z}\ .  \label{6.10kd}
\end{equation}%
Here $M$ is the Jacobian matrix of the mapping $f^{\mu }$ and $J$ is
its Jacobian.


Later on, we use the {\bf logarithmic (Newtonian) potential of
sources} $G\in L_1(\mathbb C)$ with compact supports given by the
formula:
\begin{equation} \label{eqIPOTENTIAL} {\cal N}^{G}(z)\  :=\
\frac{1}{2\pi}\int\limits_{\mathbb C} \ln|z-w|\, G(w)\ d\, m(w)\
.\end{equation}

By Lemmas 3 in \cite{GNR1} and Theorem 2 in \cite{GNR2}, we have its
basic properties:

\medskip


{\bf Proposition 5.}  {\it Let $G:\mathbb C\to\mathbb R$ have
compact support. If $G\in L_1(\mathbb C)$, then ${\cal N}^G\in
L_{r,\rm loc}(\mathbb C)$ for all $r\in [1,\infty)$, ${\cal N}^G\in
W^{1,p}_{\rm loc}(\mathbb C)$ for all $p\in [1,2)$, moreover, there
exist generalized derivatives by Sobolev $\frac{\partial^2
N^{G}}{\partial z\partial\overline z}$ and $\frac{\partial^2
N^{G}}{\partial\overline z\partial z}$ satisfying the equalities
\begin{equation}\label{eqLAP} 4\cdot\frac{\partial^2 N^{G}}{\partial
z\partial\overline z}\ =\ \triangle {\cal N}^G\ =\
4\cdot\frac{\partial^2 N^{G}}{\partial\overline z\partial z}\ =\ G\
\ \ \mbox{a.e.}
\end{equation}
Furthermore, if $G\in L_{p^{\prime}}(\mathbb C)$ for some
$p^{\prime}>1$, then ${\cal N}^G\in W^{2,p^{\prime}}_{\rm
loc}(\mathbb C)$, moreover, ${\cal N}^G\in W^{1,p}_{\rm loc}(\mathbb
C)$ for some $p>2$ and, consequently, ${\cal N}^G\in C^{\alpha}_{\rm
loc}(\mathbb C)$ with $\alpha =1-2/p$. Finally, if $G\in
L_{p^{\prime}}(\mathbb C)$ for some $p^{\prime}>2$, then ${\cal
N}^G\in C^{1,\alpha}_{\rm loc}(\mathbb C)$ with $\alpha
=1-2/p^{\prime}$.}

\medskip

As it was before, we assume here that the dilatations
$K^T_{\mu_A}(z,z_0)$ and $K_{\mu_A}(z)$ are extended by $1$ outside
of the domain $D$.

\medskip

{\bf Lemma 4.} {\it Let $D$ be a bounded domain in $\mathbb C$ with
no boundary component degenerated to a single point, $g\in
L_{p^{\prime}}(D)$, $p^{\prime}>1$, with compact support, and
$A:D\to\mathbb S^{2\times 2}$ be a matrix function with its locally
H\"older continuous entries.

\medskip

Suppose also that $K_{\mu_A}\in L^1(D)$ and, for each
$z_0\in\partial D$,
\begin{equation}\label{3omalMA}
\int\limits_{\varepsilon<|z-z_0|<\varepsilon_0}
K^T_{\mu_A}(z,z_0)\cdot\psi^2_{z_0,\varepsilon}(|z-z_0|)\,dm(z)=o(I_{z_0}^{2}(\varepsilon))\quad
\ \ \hbox{as $\varepsilon\to0$}\ ,\end{equation} where
$\varepsilon_0=\varepsilon(z_0)>0$ and $\psi_{z_0,\varepsilon}:
(0,\varepsilon_0)\to(0,\infty)$ are measurable functions
 with
\begin{equation}\label{eq3.5.3MA}
I_{z_0}(\varepsilon)\colon
=\int\limits_{\varepsilon}^{\varepsilon_0}
\psi_{z_0,\varepsilon}(t)\,dt<\infty\qquad\forall\
\varepsilon\in(0,\varepsilon_0)\,.\end{equation}

Then the Poisson type equation (\ref{eqPotential}) has the unique
weak continuous solution $u$ in the class $C^{\alpha}_{\rm loc}\cap
W^{1,p}_{\rm loc}\cap W^{2,p^{\prime}}_{\rm loc}$, $\alpha =1-2/p$,
for some $p>2$ of the Dirichlet problem (\ref{eqDIR}) in $D$ for
each continuous function $\varphi:\partial D\to{\Bbb R}$.

Moreover, $u =U\circ f$, where $f:\mathbb C\to\mathbb C$ is a $\mu_A
-$conformal mapping with $\mu_A$ extended by zero outside of $D$,
$U$ is a weak generalized harmonic function with the source $G$ of
the class $L_{p^{\prime}}(D_*)$ in the domain $D_*:=f(D)$,
\begin{equation}\label{G} G\
:=\ \frac{g}{J}\circ f^{-1}\ ,\ \ \ J(z):=|f_z|^2-|f_{\overline
z}|^2\ ,
\end{equation}
of the class $C^{\alpha}_{\rm loc}\cap W^{1,p}_{\rm loc}\cap
W^{2,p^{\prime}}_{\rm loc}$ in $D_*$, satisfying the Dirichlet
condition
\begin{equation}\label{eq_}
\lim_{w\to\zeta}\ U(w) = \varphi_*(\zeta)\ \ \  \forall\
\zeta\in\partial D_*\ ,\ \ \ \hbox{$\varphi_*:=\varphi\circ
f^{-1}|_{\partial D_*}$}\ .
\end{equation}}

Here $u$ is called a \textbf{weak solution} of the Poisson type
equation (\ref{eqPotential}) if
\begin{equation}
\int\limits_{D}\{\langle A(z)\nabla u(z),\nabla \psi \rangle \ +\
\Sigma (z)\,Q(u(z))\,\psi (z)\}\ dm_{z}\ =\ 0\ \ \ \ \forall \ \psi
\in C_{0}^{1}(D)\ .  \label{W}
\end{equation}

{\bf Remark 14.} In turn, by the proof below, $U := {\cal H} + {\cal
N}^G$, where ${\cal H}$ is the unique harmonic function in $D_*$,
satisfying the Dirichlet condition
\begin{equation}\label{eq*}
\lim_{w\to\zeta}\ {\cal H}(w) = \varphi^*(\zeta)\ \ \  \forall\
\zeta\in\partial D_*\ ,\ \ \ \hbox{$\varphi^*:=\varphi_* - {\cal
N}^G|_{\partial D_*}$}\ .
\end{equation}
Furthermore, arguing similarly to the proof, we obtain also by
Proposition 5 that if $g\in L_{p^{\prime}}(D)$ for some
$p^{\prime}>2$, then in addition $u\in C^{1,\alpha^{\prime}}_{\rm
loc}(D)$ with $\alpha^{\prime} =1-2/p^{\prime}$. In the case, the
function $U$ is a generalized harmonic function with the source
$G\in L_{p^{\prime}}(D_*)$ of the class $C^{1,\alpha^{\prime}}_{\rm
loc}\cap W^{2,p^{\prime}}_{\rm loc}$ in the domain $D_*$.

\begin{proof} Let $f:\mathbb{C}\rightarrow {\mathbb{C}}$ be a $\mu_A $-conformal
mapping from Proposition 1 with the complex coefficient $\mu_A$ in
$\mathbb{C}$ extended by zero outside of $D$. Since entries of $A$
are locally H\"older continuous, the mapping $f|_D$ is smooth, see
e.g. \cite{Iw} and \cite{IwDis}, and, moreover, see e.g. Theorem
V.7.1 in \cite{LV}, its continuous Jacobian
\begin{equation}
J(z)\ =\ |f_{z}|^{2}\,-\,|f_{\bar{z}}|^{2}\ >\ 0\ \ \ \ \ \ \
\forall \ z\in D\ .  \label{6.13_2}
\end{equation}Consequently, $f^{-1}$
is also smooth in the domain $D_*:=f(D)$, see e.g. formulas I.C(3)
in \cite{Alf}. Thus, by the replacement of variables, see e.g. the
point (vi) of Theorem 5 in \cite{ABe}, the function $G$ in (\ref{G})
belongs to the class $L_{p^{\prime}}(D_*)$ because of $G$ has
compact support in $D_*$ by hypotheses of the lemma and in view of
homeomorphism of $f$, and because of the continuous function
$J^{-1}\circ f^{-1}$ is bounded over the support of $G$.

Next, the domain $D_*$ is bounded and has no boundary component
de\-ge\-ne\-ra\-ted to a single point because of $D$ is so by
hypotheses of the lemma and because of the mapping $f$ is a
homeomorphism of $\mathbb C$ into itself. Thus, by Corollary 4.1.8
and Theorem 4.2.2 in \cite{Rans}, there is the unique harmonic
function ${\cal H}:D\to\mathbb R$, satisfying the Dirichlet
condition (\ref{eq*}). Thus, by Proposition 5 $U:={\cal H}+{\cal
N}^G$ is a weak generalized harmonic function with the source $G$ of
the class $L_{p^{\prime}}(D_*)$ in the domain $D_*$, satisfying the
Dirichlet condition (\ref{eq_}). Note that again by Proposition 5
$U\in W^{2,p^{\prime}}_{\rm loc}(\mathbb C)$, moreover, $U\in
W^{1,p}_{\rm loc}(\mathbb C)$ for some $p>2$ and, consequently,
$U\in C^{\alpha}_{\rm loc}(\mathbb C)$ with $\alpha =1-2/p$.

Finally, by Proposition 3.1 in \cite{GNR}, see (\ref{6.10kd}), the
function $u:=U\circ f$ gives the desired solution of the Poisson
type equation (\ref{eqPotential}) because $f|_D$ is a local
quasi-isometry in $D$ of the class $C^{1}$, see e.g. 1.1.7 in
\cite{Maz}, and such a solution is unique.
\end{proof}


{\bf Remark 15.} Note that if the family of the functions
$\psi_{z_0,\varepsilon}(t)\equiv\psi_{z_0}(t)$ is independent on the
parameter $\varepsilon$, then the condition (\ref{3omalMA}) implies
that $I_{z_0}(\varepsilon)\to \infty$ as $\varepsilon\to 0$. This
follows immediately from arguments by contradiction, apply for it
(\ref{eqConnect}) and the condition $K_{\mu_A}\in L^1(D)$. Note also
that (\ref{3omalMA}) holds, in particular, if, for some
$\varepsilon_0=\varepsilon(z_0)$,
\begin{equation}\label{333omalMA}
\int\limits_{|z-z_0|<\varepsilon_0}
K^T_{\mu_A}(z,z_0)\cdot\psi_{z_0}^2(|z-z_0|)\,dm(z)<\infty \qquad
\forall\ z_0\in\partial{D}\end{equation} and
$I_{z_0}(\varepsilon)\to \infty$ as $\varepsilon\to 0$. In other
words, for the existence of regular enough weak solutions of the
Dirichlet problem (\ref{eqDIR}) in $D$ to the Poisson type equation
(\ref{eqPotential}) with arbitrary continuous boundary functions
$\varphi$, it is sufficient that the integral in (\ref{333omalMA})
converges for some nonnegative function $\psi_{z_0}(t)$ that is
locally integrable over $(0,\varepsilon_0 ]$ but has a nonintegrable
singularity at $0$. The functions $\log^{\lambda}(e/|z-z_0|)$,
$\lambda\in (0,1)$, $z\in\Bbb D$, $z_0\in\overline{\Bbb D}$, and
$\psi(t)=1/(t \,\, \log(e/t))$, $t\in(0,1)$, show that the condition
(\ref{333omalMA}) is compatible with the condition
$I_{z_0}(\varepsilon)\to\infty $ as $\varepsilon\to 0$. Furthermore,
the condition (\ref{3omalMA}) in Lemma 4 shows that, for the
existence of such solutions of the Dirichlet problem (\ref{eqDIR})
to the Poisson type equation (\ref{eqPotential}), it is sufficient
even that the integral in (\ref{333omalMA}) to be divergent in a
controlled way.

\medskip

Similarly to Section 3, we derive from Lemma 4 the next series of
results.

\medskip

{\bf Theorem 9.} {\it Let $D$ be a bounded domain in $\mathbb C$
with no boundary component degenerated to a single point, $g\in
L_{p^{\prime}}(D)$, $p^{\prime}>1$, with compact support, and
$A:D\to\mathbb S^{2\times 2}$ be a matrix function whose elements
$a_{ij}(z)$, $i,j=1,2$ are measurable and locally bounded,
$K_{\mu_A}\in L^1(D)$ and, for each $z_0\in \partial{D}$,
$K^T_{\mu_A}(z,z_0)\leqslant Q_{z_0}(z)$ in its neighborhood
$U_{z_0}$ for a function $Q_{z_0}: U_{z_0}\to[0,\infty]$ in the
class ${\rm FMO}({z_0})$.

Then the Poisson type equation (\ref{eqPotential}) has the unique
weak continuous solution $u$ in the class $C^{\alpha}_{\rm loc}\cap
W^{1,p}_{\rm loc}\cap W^{2,p^{\prime}}_{\rm loc}$, $\alpha =1-2/p$,
for some $p>2$ of the Dirichlet problem (\ref{eqDIR}) in $D$ for
each continuous function $\varphi:\partial D\to{\Bbb R}$.

Moreover, $u =U\circ f$, $U := {\cal H} + {\cal N}^G$, where
$f:\mathbb C\to\mathbb C$ is a $\mu_A -$conformal mapping with
$\mu_A$ extended by zero outside of $D$, $U$ is a weak generalized
harmonic function in $D_*:=f(D)$ with the source $G\in
L_{p^{\prime}}(D_*)$ calculated in (\ref{G}) and ${\cal H}$ is the
unique harmonic function with the Dirichlet condition (\ref{eq*}).}

\medskip

{\bf Corollary 17.} {\it Let $D$ be a bounded domain in $\mathbb C$
with no boundary component degenerated to a single point, $g\in
L_{p^{\prime}}(D)$, $p^{\prime}>1$, with compact support, and
$A:D\to\mathbb S^{2\times 2}$ be a matrix function whose elements
$a_{ij}(z)$, $i,j=1,2$ are measurable and locally bounded,
$K_{\mu_A}\in L^1(D)$ and, for each $z_0\in \partial{D}$,
\begin{equation}\label{eqMEANmA}\overline{\lim\limits_{\varepsilon\to0}}\quad
\dashint_{B(z_0,\varepsilon)}K^T_{\mu_A}(z,z_0)\,dm(z)<\infty\
.\end{equation}

Then all the conclusions of Theorem 9 on solutions of the Dirichlet
problem (\ref{eqDIR}) with continuous data $\varphi:\partial
D\to{\Bbb R}$ to the Poisson type equation (\ref{eqPotential})
hold.}

\medskip

{\bf Corollary 18.} {\it Let $D$ be a bounded domain in $\mathbb C$
with no boundary component degenerated to a single point, $g\in
L_{p^{\prime}}(D)$, $p^{\prime}>1$, with compact support, and
$A:D\to\mathbb S^{2\times 2}$ be a matrix function whose elements
$a_{ij}(z)$, $i,j=1,2$ are measurable and locally bounded in $D$,
and $K_{\mu_A}$ have a dominant $Q:U\to[1,\infty)$ of the class {\rm
BMO}$_{\rm loc}(U)$ in a neighborhood $U$ of $\partial D$. Then all
the conclusions of Theorem 9 hold.}

\medskip

{\bf Remark 16.} In particular, the conclusions of Theorem 9 hold if
$Q\in{\rm W}^{1,2}_{\rm loc}$.

\medskip

{\bf Corollary 19.} {\it Let $D$ be a bounded domain in $\mathbb C$
with no boundary component degenerated to a single point, $g\in
L_{p^{\prime}}(D)$, $p^{\prime}>1$, with compact support, and
$A:D\to\mathbb S^{2\times 2}$ be a matrix function whose elements
$a_{ij}(z)$, $i,j=1,2$ are measurable and locally bounded in $D$,
and $K_{\mu_A}$ have a dominant $Q:U\to[1,\infty)$ of the class {\rm
FMO}$(U)$ in a neighborhood $U$ of $\partial D$. Then all the
conclusions of Theorem 9 hold. }

\medskip

{\bf Theorem 10.} {\it Let $D$ be a bounded domain in $\mathbb C$
with no boundary component degenerated to a single point, $g\in
L_{p^{\prime}}(D)$, $p^{\prime}>1$, with compact support, and
$A:D\to\mathbb S^{2\times 2}$ be a matrix function whose elements
$a_{ij}(z)$, $i,j=1,2$ are measurable and locally bounded in $D$,
$K_{\mu_A}\in L^1(D)$, and, for each $z_0\in\partial{D}$,
$\varepsilon_0=\varepsilon(z_0)>0$,
\begin{equation}\label{eqLOGmA}
\int\limits_{\varepsilon<|z-z_0|<\varepsilon_0}K^T_{\mu_A}(z,z_0)\,\frac{dm(z)}{|z-z_0|^2}
=o\left(\left[\log\frac{1}{\varepsilon}\right]^2\right)\qquad\hbox{as
$\varepsilon\to 0$}\ .\end{equation}

Then the Poisson type equation (\ref{eqPotential}) has the unique
weak continuous solution $u$ in the class $C^{\alpha}_{\rm loc}\cap
W^{1,p}_{\rm loc}\cap W^{2,p^{\prime}}_{\rm loc}$, $\alpha =1-2/p$,
for some $p>2$ of the Dirichlet problem (\ref{eqDIR}) in $D$ for
each continuous function $\varphi:\partial D\to{\Bbb R}$.

Moreover, $u =U\circ f$, $U := {\cal H} + {\cal N}^G$, where
$f:\mathbb C\to\mathbb C$ is a $\mu_A -$conformal mapping with
$\mu_A$ extended by zero outside of $D$, $U$ is a weak generalized
harmonic function in $D_*:=f(D)$ with the source $G\in
L_{p^{\prime}}(D_*)$ calculated in (\ref{G}) and ${\cal H}$ is the
unique harmonic function with the Dirichlet condition (\ref{eq*}).}

\medskip

{\bf Remark 17.} Choosing in Lemma 4 the function
$\psi(t)=1/(t\log{1/t})$ instead of $\psi(t)=1/t$, we are able to
replace (\ref{eqLOGmA}) by
\begin{equation}\label{eqLOGLOGmA}
\int\limits_{\varepsilon<|z-z_0|<\varepsilon_0}\frac{K^T_{\mu_A}(z,z_0)\,dm(z)}
{\left(|z-z_0|\log{\frac{1}{|z-z_0|}}\right)^2}
=o\left(\left[\log\log\frac{1}{\varepsilon}\right]^2\right)\end{equation}
In general, we are able to give here the whole scale of the
corresponding conditions using functions $\psi(t)$ of the form
$1/(t\log{1}/{t}\cdot\log\log{1}/{t}\cdot\ldots\cdot\log\ldots\log{1}/{t})$.

\medskip

{\bf Theorem 11.} {\it Let $D$ be a bounded domain in $\mathbb C$
with no boundary component degenerated to a single point, $g\in
L_{p^{\prime}}(D)$, $p^{\prime}>1$, with compact support, and
$A:D\to\mathbb S^{2\times 2}$ be a matrix function whose elements
$a_{ij}(z)$, $i,j=1,2$ are measurable and locally bounded in $D$,
$K_{\mu_A}\in L^1(D)$, and, for each $z_0\in\partial{D}$,
$\varepsilon_0=\varepsilon(z_0)>0$,
\begin{equation}\label{eqLEHTOmA}\int\limits_{0}^{\varepsilon_0}
\frac{dr}{rk^T_{\mu_A}(z_0,r)}\ =\ \infty\ ,\end{equation} where
$k_{\mu_A}^T(z_0,r)$ is the integral mean of $K^T_{{\mu_A}}(z,z_0)$
over the circle $S(z_0,r)$.

\medskip

Then the Poisson type equation (\ref{eqPotential}) has the unique
weak continuous solution $u$ in the class $C^{\alpha}_{\rm loc}\cap
W^{1,p}_{\rm loc}\cap W^{2,p^{\prime}}_{\rm loc}$, $\alpha =1-2/p$,
for some $p>2$ of the Dirichlet problem (\ref{eqDIR}) in $D$ for
each continuous function $\varphi:\partial D\to{\Bbb R}$.

Moreover, $u =U\circ f$, $U := {\cal H} + {\cal N}^G$, where
$f:\mathbb C\to\mathbb C$ is a $\mu_A -$conformal mapping with
$\mu_A$ extended by zero outside of $D$, $U$ is a weak generalized
harmonic function in $D_*:=f(D)$ with the source $G\in
L_{p^{\prime}}(D_*)$ calculated in (\ref{G}) and ${\cal H}$ is the
unique harmonic function with the Dirichlet condition (\ref{eq*}).}

\medskip

{\bf Corollary 20.} {\it Let $D$ be a bounded domain in $\mathbb C$
with no boundary component degenerated to a single point, $g\in
L_{p^{\prime}}(D)$, $p^{\prime}>1$, with compact support, and
$A:D\to\mathbb S^{2\times 2}$ be a matrix function whose elements
$a_{ij}(z)$, $i,j=1,2$ are measurable and locally bounded in $D$,
$K_{\mu_A}\in L^1(D)$, and, for each $z_0\in\partial{D}$,
\begin{equation}\label{eqLOGkMA}k^T_{\mu_A}(z_0,\varepsilon)=O\left(\log\frac{1}{\varepsilon}\right)
\qquad\mbox{as}\ \varepsilon\to0\ .\end{equation}

Then all the conclusions of Theorem 11 on solutions of the Dirichlet
problem (\ref{eqDIR}) with continuous data $\varphi:\partial
D\to{\Bbb R}$ to the Poisson type equation (\ref{eqPotential})
hold.}

\medskip

{\bf Remark 18.} In particular, all the conclusions of Theorem 11
hold if, for each point $z_0\in\partial{D}$,
\begin{equation}\label{eqLOGKmA} K^T_{\mu_A}(z,z_0)=O\left(\log\frac{1}{|z-z_0|}\right)\qquad{\rm
as}\quad z\to z_0\ .\end{equation} Moreover, (\ref{eqLOGkMA}) can be
replaced by the whole series of weaker condition
\begin{equation}\label{edLOGLOGkMA}
k^T_{\mu_A}(z_0,\varepsilon)=O\left(\left[\log\frac{1}{\varepsilon}\cdot\log\log\frac{1}
{\varepsilon}\cdot\ldots\cdot\log\ldots\log\frac{1}{\varepsilon}
\right]\right) \qquad\forall\ z_0\in \partial{D}\ .
\end{equation}

\medskip

{\bf Theorem 12.} {\it Let $D$ be a bounded domain in $\mathbb C$
with no boundary component degenerated to a single point, $g\in
L_{p^{\prime}}(D)$, $p^{\prime}>1$, with compact support, and
$A:D\to\mathbb S^{2\times 2}$ be a matrix function whose elements
$a_{ij}(z)$, $i,j=1,2$ are measurable and locally bounded in $D$,
$K_{\mu_A}\in L^1(D)$, and, for each $z_0\in\partial{D}$ and  a
neighborhood $U_{z_0}$ of $z_0$,
\begin{equation}\label{eqINTEGRALmA}\int\limits_{U_{z_0}}\Phi_{z_0}\left(K^T_{\mu_A}(z,z_0)\right)\,dm(z)\ <\ \infty\
,
\end{equation} where  $\Phi_{z_0}:[0,\infty]\to[0,\infty]$ is a convex
non-decreasing function such that, for some $\Delta(z_0)>0$,
\begin{equation}\label{eqINTmA}
\int\limits_{\Delta(z_0)}^{\infty}\log\,\Phi_{z_0}(t)\,\frac{dt}{t^2}\
=\ +\infty\ .\end{equation}

Then the Poisson type equation (\ref{eqPotential}) has the unique
weak continuous solution $u$ in the class $C^{\alpha}_{\rm loc}\cap
W^{1,p}_{\rm loc}\cap W^{2,p^{\prime}}_{\rm loc}$, $\alpha =1-2/p$,
for some $p>2$ of the Dirichlet problem (\ref{eqDIR}) in $D$ for
each continuous function $\varphi:\partial D\to{\Bbb R}$.

Moreover, $u =U\circ f$, $U := {\cal H} + {\cal N}^G$, where
$f:\mathbb C\to\mathbb C$ is a $\mu_A -$conformal mapping with
$\mu_A$ extended by zero outside of $D$, $U$ is a weak generalized
harmonic function in $D_*:=f(D)$ with the source $G\in
L_{p^{\prime}}(D_*)$ calculated in (\ref{G}) and ${\cal H}$ is the
unique harmonic function with the Dirichlet condition (\ref{eq*}).}

\bigskip

{\bf Corollary 21.} {\it Let $D$ be a bounded domain in $\mathbb C$
with no boundary component degenerated to a single point, $g\in
L_{p^{\prime}}(D)$, $p^{\prime}>1$, with compact support, and
$A:D\to\mathbb S^{2\times 2}$ be a matrix function whose elements
$a_{ij}(z)$, $i,j=1,2$ are measurable and locally bounded in $D$,
$K_{\mu_A}\in L^1(D)$, and, for each $z_0\in\partial{D}$,  a
neighborhood $U_{z_0}$ of $z_0$ and $\alpha(z_0)>0$,
\begin{equation}\label{eqEXPmA}\int\limits_{U_{z_0}}e^{\alpha(z_0) K^T_{\mu_A}(z,z_0)}\,dm(z)<\infty
\ .\end{equation}

Then all the conclusions of Theorem 12 on solutions of the Dirichlet
problem (\ref{eqDIR}) with continuous data $\varphi:\partial
D\to{\Bbb R}$ to the Poisson type equation (\ref{eqPotential})
hold.}

{\bf Corollary 22.} {\it Let $D$ be a bounded domain in $\mathbb C$
with no boundary component degenerated to a single point, $g\in
L_{p^{\prime}}(D)$, $p^{\prime}>1$, with compact support, and
$A:D\to\mathbb S^{2\times 2}$ be a matrix function whose elements
$a_{ij}(z)$, $i,j=1,2$ are measurable and locally bounded in $D$,
$K_{\mu_A}\in L^1(D)$, and, for a neighborhood $U$ of $\partial D$,
\begin{equation}\label{eqINTKmA}\int\limits_{U}\Phi\left(K_{\mu_A}(z)\right)\,dm(z)<\infty\ ,\end{equation}
where $\Phi:[0,\infty]\to[0,\infty]$ is a convex non-decreasing
function such that, for some $\delta>0$,
\begin{equation}\label{eqINTFmA}
\int\limits_{\delta}^{\infty}\log\,\Phi(t)\,\frac{dt}{t^2}\ =\
+\infty\ .\end{equation}

Then all the conclusions of Theorem 12 on solutions of the Dirichlet
problem (\ref{eqDIR}) with continuous data $\varphi:\partial
D\to{\Bbb R}$ to the Poisson type equation (\ref{eqPotential})
hold.}

\medskip

{\bf Remark 19.} By Theorems 2.5 and 5.1 in \cite{RSY$_5$},
condition (\ref{eqINTFmA}) is not only sufficient but also necessary
to have a regular enough weak solution $u$ of the Dirichlet problem
(\ref{eqDIR}) in $D$ for all the Poisson type equations
(\ref{eqPotential}), satisfying the integral constraints
(\ref{eqINTKmA}), for arbitrary continuous functions
$\varphi:\partial D\to\Bbb{R}$ because such solutions have the
representation through regular homeomorphic solutions $f=f^{\mu}$ of
the homogeneous Beltrami equation (\ref{a}) with $\mu=\mu_A$.

\medskip

{\bf Corollary 23.} {\it Let $D$ be a bounded domain in $\mathbb C$
with no boundary component degenerated to a single point, $g\in
L_{p^{\prime}}(D)$, $p^{\prime}>1$, with compact support, and
$A:D\to\mathbb S^{2\times 2}$ be a matrix function whose elements
$a_{ij}(z)$, $i,j=1,2$ are measurable and locally bounded in $D$,
$K_{\mu_A}\in L^1(D)$ and, for a neighborhood $U$ of $\partial{D}$
and $\alpha>0$,
\begin{equation}\label{eqEXPAmA}\int\limits_{D}e^{\alpha K_{\mu_A}(z)}\,dm(z)\ <\
\infty\ .
\end{equation}

Then all the conclusions of Theorem 12 on solutions of the Dirichlet
problem (\ref{eqDIR}) with continuous data $\varphi:\partial
D\to{\Bbb R}$ to the Poisson type equation (\ref{eqPotential})
hold.}

\medskip

As a result, we have a number of effective integral criteria for the
solvability of the classical Dirichlet problem (\ref{eqDIR}) in the
most general admissible domains to one of the main equations
(\ref{eqPotential}) of the hydromechanics (fluid mechanics) in
anisotropic and inhomogeneous media.

\bigskip

{\bf Acknowledgments.} The first 3 authors are partially supported
by the project "Mathematical modelling of complex dynamical systems
and processes caused by the state security", No. 0123U100853, of
National Academy of Scien\-ces of Ukraine and by the Grant
EFDS-FL2-08 of the found of the European Federation of Academies of
Sciences and Humanities (ALLEA).


\bigskip

{\bf \noindent Vladimir Gutlyanskii, Olga Nesmelova, Vladimir Ryazanov} \\
Institute of Applied Mathematics and Mechanics of NASU, Slavyansk,\\
Institute of Mathematics of National Academy of Sciences of Ukraine, Kiev,\\
vgutlyanskii@gmail.com, star-o@ukr.net, vl.ryazanov1@gmail.com

\bigskip

\noindent {\bf Eduard Yakubov}\\
Holon Institute of Technology, Holon, Israel,\\
yakubov@hit.ac.il, eduardyakubov@gmail.com


\begin{thebibliography}{100}

\bibitem{Alf}
{Ahlfors L.} (1966) \emph{Lectures on Quasiconformal Mappings.} New
York: Van Nostrand.

\bibitem{ABe}
Ahlfors, L.V., Bers, L. (1960) Riemann’s mapping theorem for
variable metrics. {\it Ann. Math. (2) 72}, 385-404.

\bibitem{And}
{\sc Andreian Cazacu C.:} On the length-area dilatation. - Complex
Variables, Theory Appl. 50:7--11, 2005, 765--776.

\bibitem{AIM}
{Astala K., Iwaniec T., Martin G.J.} (2009) \emph{Elliptic
differential equations and quasiconformal mappings in the plane.}
Princeton Math. Ser.  {\bf 48}. Princeton: Princeton Univ. Press.

\bibitem{Bojar} {\sc
Bojarski B.:} Generalized solutions of a system of differential
equations of the first order of the elliptic type with discontinuous
coefficients. - Mat. Sb., N. Ser. 43(85):4, 1958, 451–503; transl.
in Report of Univ. of Jyv\"askyl\"a, Dept. Math. and Stat., 118. -
Univ. of Jyv\"askyl\"a, 2009.

\bibitem{BGMR} {\sc Bojarski B., Gutlyanskii V., Martio O., Ryazanov V.:}
In\-fi\-ni\-te\-si\-mal geometry of quasiconformal and bi-Lipschitz
mappings in the plane. - EMS Tracts in Mathematics, 19.  Z\"urich:
European Mathematical Society (EMS), 2013.

\bibitem{BGR$_1$} {\sc Bojarski B., Gutlyanskii V., Ryazanov V.:}
On integral conditions for the general Beltrami equations. - Complex
Anal. Oper. Theory 5:3, 2011, 835-845.

\bibitem{BGR$_2$} {\sc Bojarski B., Gutlyanskii V., Ryazanov V.:}
On existence and representation of solutions for general degenerate
Beltrami equations. - Complex Var. Elliptic Equ. 59:1, 2014, 67-75.

\bibitem{BN}
{\sc Brezis H., Nirenberg L.:} {Degree theory and BMO. I. Compact
manifolds without boundaries.} Selecta Math. (N.S.) 1:2, 1995,
197--263.

\bibitem{CFL}
{\sc Chiarenza F., Frasca M., Longo P.:} {$W^{2,p}$-solvability of
the Dirichlet problem for nondivergence elliptic equations with VMO
coefficients.} - Trans. Amer. Math. Soc. 336:2, 1993, 841--853.

\bibitem{GMSV} {\sc
Gutlyanskii V., Martio O., Sugawa T., Vuorinen M.:} On the
degenerate Beltrami equation. - Trans. Amer. Math. Soc. 357:3, 2005,
875-900.

\bibitem{GNR-}
V. Gutlyanskii, O. Nesmelova, and V. Ryazanov, ``On a
model semilinear elliptic equation in the plane,'' \emph{J. Math. Sci.,} {\bf 220}%
(5), 603--614 (2017); transl. from \emph{Ukr. Mat. Visn.,} {\bf
13}(1), 91--105 (2016).

\bibitem{GNR}
V. Gutlyanskii, O.  Nesmelova, and V. Ryazanov, ``On quasiconformal
maps and semi-linear equations in the plane,'' \emph{J. Math. Sci.,}
{\bf 229}(1), 7--29 (2018); transl. from. \emph{Ukr. Mat. Visn.,}
{\bf 14}(2), 161--191 (2017)

\bibitem{GNR1}
Gutlyanskii, V., Nesmelova, O., Ryazanov, V. (2019). To the theo\-ry
of semi-linear equations in the plane. {\it J. Math. Sci.} (USA)
242(6), 833--859; transl. (2019) from {\it Ukr. Mat. Visn.}, 16, No.
1, 105--140.

\bibitem{GNR2}
Gutlyanskii, V., Nesmelova, O., Ryazanov, V.  (2020). On a
quasilinear Poisson equation in the plane. {\it Anal. Math. Phys.,
10}(1), Paper No. 6, 1--14.

\bibitem{GNRY*}
Gutlyanskii, V., Nesmelova, O., Ryazanov, V., Yakubov, E. (2023). On
the Hilbert problem for semi-linear Beltrami equations. {\it J.
Math. Sci. (USA) 270}, no. 3, 428-448; transl. from {\it Ukr. Mat.
Visn. 19}, No. 4, 489 – 516; {see also arXiv: 2206.05045 [math.AP]}
https://doi.org/10.48550/arXiv.2206.05045






\bibitem{GNRY}
Gutlyanskii, V., Nesmelova, O., Ryazanov, V., Yefimushkin, A.
(2021). Logarithmic potential and generalized analytic functions.
{\it J. Math. Sci. (USA) 256}, No. 6, 735-752; transl. from {\it
Ukr. Mat. Visn. 18}, No. 1, 12-36.

\bibitem{GRSY*} {\sc Gutlyanskii V., Ryazanov V., Sevostyanov E., Yakubov E.}
(2022). BMO and Dirichlet problem for degenerate Beltrami equation.
J. Math. Sci. (USA) 268, no. 2, 157–177; transl. from {\it Ukr. Mat.
Visn. 19}, No. 3, 327–354.

\bibitem{GRSY} {\sc Gutlyanskii V., Ryazanov V., Srebro U., Yakubov E.:}
The Beltrami Equation: A Geometric Approach. - Developments in
Mathematics 26, Springer: Berlin, 2012.

\bibitem{HKM}
Heinonen, J., Kilpel\"ainen, T., Martio, O. (1993). \emph{Nonlinear
potential theory of degenerate elliptic equations}. Oxford
Mathematical Monographs. Oxford Science Publications, The Clarendon
Press, Oxford University Press, New York.

\bibitem{IR}
{\sc Ignat'ev A.A., Ryazanov V.I.:} {Finite mean oscillation in the
mapping theory.} - Ukrainian Math. Bull. 2:3, 2005, 403-424.

\bibitem{Iw}
T. Iwaniec,  ``Regularity of solutions of certain degenerate
elliptic systems of equations that realize quasiconformal mappings
in n-dimensional space,'' \emph{Differential and integral equations.
Boundary value problems}, 97--111 (1979).

\bibitem{IwDis}
T. Iwaniec, ``Regularity theorems for solutions of partial
differential equations for quasiconformal mappings in several
dimensions,'' \emph{Dissertationes Math. (Rozprawy Mat.),} {\bf
198}, 45 pp. (1982).

\bibitem{IS}
{\sc Iwaniec T., C. Sbordone C.:} {Riesz transforms and elliptic
PDEs with VMO coefficients.} - J. Anal. Math. 74, 1998, 183--212.

\bibitem{JN}
{\sc John F., Nirenberg L.:} {On functions of bounded mean
oscillation.} - Comm. Pure Appl. Math. 14, 1961, 415--426.

\bibitem{Ko}
Koosis, P. (1998). \emph{Introduction to $H^p$ spaces.} Cambridge
Tracts in Ma\-the\-ma\-tics. 115. Cambridge Univ. Press, Cambridge.

\bibitem{Le}
{\sc Lehto O.:} Homeomorphisms with a prescribed dilatation. -
Lecture Notes in Math. 118, 1968, 58-73.

\bibitem{LV}
{Lehto O., Virtanen K.J.} (1973) \emph{Quasiconformal mappings in
the plane.} Springer-Verlag: Berlin, Heidelberg.

\bibitem{MRV}
{\sc Martio O., Ryazanov V., Vuorinen M.:} {BMO and Injectivity of
Space Quasiregular Mappings.} - Math. Nachr. 205, 1999, 149--161.

\bibitem{MRSY} {\sc Martio O., Ryazanov V., Srebro U., Yakubov E.:} Moduli in modern mapping
theory. - Springer Monographs in Mathematics. - Springer: New York,
2009.

\bibitem{Maz}
V. G. Maz'ja,  \emph{Sobolev spaces}, Springer-Verlag, Berlin, 1985.

\bibitem{Ne} {\sc Nevanlinna R.:} Eindeutige analytische Funktionen. 2. Aufl. Reprint.
(German) - Die Grundlehren der mathematischen Wissenschaften. Band
46. Springer-Verlag: Berlin-Heidelberg-New York, 1974.

\bibitem{Pal}
{\sc Palagachev D.K.:} {Quasilinear elliptic equations with VMO
coefficients.} - Trans. Amer. Math. Soc. 347:7, 1995, 2481--2493.

\bibitem{Ra$_1$}
{\sc Ragusa M.A.:} {Elliptic boundary value problem in vanishing
mean oscillation hypothesis.} - Comment. Math. Univ. Carolin. 40:4,
1999, 651--663.

\bibitem{Ra$_2$}
{\sc Ragusa M.A., Tachikawa A.:} Partial regularity of the
minimizers of quadratic functionals with VMO coefficients. J. Lond.
Math. Soc., II. Ser. 72:3, 2005, 609-620.

\bibitem{Rans} {\sc Ransford Th.:} Potential theory in the complex plane.
- London Mathematical Society Student Texts
28, Univ. Press: Cambridge, 1995.

\bibitem{RR}
{\sc Reimann H.M., Rychener T.:} {Funktionen Beschr\"ankter
Mittlerer Oscillation.} - Lecture Notes in Math. 487, 1975.


\bibitem{R4}
Ryazanov, V. (2021). On Hilbert and Riemann problems for
ge\-ne\-ra\-li\-zed ana\-ly\-tic functions and applications.
\textit{Anal. Math. Phys., 11}(5). Published online: 22 November
2020.


\bibitem{RS}
Ryazanov, Vladimir I.; Salimov, Ruslan R. Weakly planar spaces and
boundaries in the theory of mappings. (Russian) Ukr. Mat. Visn. 4
(2007), no. 2, 199–234, 307; translation in Ukr. Math. Bull. 4
(2007), no. 2, 199–234

\bibitem{RSY} {\sc Ryazanov V., Srebro U., Yakubov E.:} BMO-quasiconformal mappings.
- J. d'Anal. Math. 83, 2001, 1-20.

\bibitem{RSY$_2$}
{\sc Ryazanov V., Srebro U., Yakubov E.:} {Beltrami equation and FMO
functions.} - Contemp. Math. 382, Israel Math. Conf. Proc., 2005,
357--364.

\bibitem{RSY$_3$} {\sc Ryazanov~V., U. Srebro, and E. Yakubov:}
Finite mean oscillation and the Beltrami equation. - Israel Math. J.
153, 2006, 247--266.

\bibitem{RSY$_4$} {\sc Ryazanov~V., U. Srebro, and E. Yakubov:}
On the theory of the Beltrami equation. - Ukr. Math. J. 58:11, 2006,
1786-1798.

\bibitem{RSY$_5$} {\sc Ryazanov V., Srebro U., Yakubov E.:}
Integral conditions in the theory of the Beltrami equations. Complex
Var. Elliptic Equ. 57:12, 2012, 1247-1270.


\bibitem{ST}
{\sc Saff E.B., Totik V.:} Logarithmic potentials with external
fields. - Grundlehren der Mathematischen Wissenschaften. 316.
Springer: Berlin, 1997.

\bibitem{Sar}
{\sc Sarason D.:} {Functions of vanishing mean oscillation.} -
Trans. Amer. Math. Soc. 207, 1975, 391--405.

\bibitem{So}
Sobolev, S.L. (1963). \emph{Applications of functional analysis in
mathematical physics.} Transl. of Math. Mon. 7. AMS, Providence,
R.I.

\bibitem{Vek}
Vekua, I.N. (1962). \emph{Generalized analytic functions.} Pergamon
Press. London-Paris-Frankfurt; Addison-Wesley Publishing Co., Inc.,
Reading, Mass.

\bibitem{Wien}
{\sc Wiener N.} The Dirichlet problem. Mass. J. of Math. 3, 1924,
129-146.

\end{thebibliography}
\end{document}